\documentclass{article}

\usepackage{etoolbox} 
\usepackage[shortlabels,inline]{enumitem}

\usepackage[immediate]{silence}
\usepackage{environ}
\usepackage{xparse}
\usepackage{listofitems}

\usepackage[
    a4paper,
    margin=20mm 
]{geometry}
\usepackage[onehalfspacing]{setspace} 

\usepackage{sectsty} 
\DeactivateWarningFilters[temp] 

\usepackage{titlesec}

\usepackage{amsmath, amssymb, amsthm}
\usepackage{thmtools, thm-restate}
\usepackage{mathtools}
\usepackage[a]{esvect}

\usepackage{array, makecell}

\usepackage{graphicx}
\usepackage[dvipsnames]{xcolor}
\usepackage[pstarrows]{pict2e}

\usepackage[
    setpagesize=false,
    colorlinks=true,
    linkcolor=blue,
    citecolor=blue,
    pdfencoding=auto,
    psdextra,
]{hyperref}
\usepackage{cleveref}
\usepackage[
    sorting=nyt,
    maxbibnames=99,
    url=false 
]{biblatex}

\usepackage[
    deletedmarkup=sout,
    commentmarkup=uwave,
    authormarkuptext=name,
]{changes}
\usepackage[linewidth=0.4pt]{mdframed}

\AtEveryBibitem{%
    \clearfield{doi}
    \clearfield{issn}
    \clearfield{eprintclass} 
}

\setlist[enumerate,1]{label=(\arabic*), ref=(\arabic*)}
\setlist[enumerate,3]{label=(\roman*), ref=(\roman*)}

\allsectionsfont{\boldmath} 

\titlespacing*{\paragraph}{0pt}{1.5ex plus 1ex minus .2ex}{1em}
\titlespacing*{\subparagraph}{0pt}{1.5ex plus 1ex minus .2ex}{1em}

\newcounter{propcounter}
\newenvironment{proplist}[1][]{%
    \stepcounter{propcounter}%
    \enumerate[label = {\bfseries \Alph{propcounter}\arabic{enumi}}, #1]%
}
{\endenumerate}

\crefname{subsection}{subsection}{subsections}

\makeatletter
\newcommand{\newreptheorem}[2]{
    \newtheorem*{rep@#1}{\rep@title}%
    \newenvironment{rep#1}[1]{%
            \def\rep@title{#2 \ref*{##1}}%
            \begin{rep@#1}%
        }
        {\end{rep@#1}}%
}
\makeatother

\theoremstyle{plain}
\newtheorem{theorem}{Theorem}[section]
\newtheorem{lemma}[theorem]{Lemma}
\newtheorem{corollary}[theorem]{Corollary}
\newtheorem{proposition}[theorem]{Proposition}

\newtheorem{question}[theorem]{Question}

\newtheorem{claim}{Claim}[theorem]
\newtheorem*{claim*}{Claim}

\theoremstyle{definition}
\newtheorem{definition}[theorem]{Definition}

\newtheorem*{remark*}{Remark}

\makeatletter
\newenvironment{claimproof}[1][Proof]{\par
	\pushQED{\qed}%
	
	\normalfont \topsep6\p@\@plus6\p@\relax
	\trivlist
	\item[\hskip\labelsep
	\textit{#1}\@addpunct{.}~]\ignorespaces
}{%
	\popQED\endtrivlist\@endpefalse
}
\makeatother

\newlist{Cases}{enumerate}{3}
\setlist[Cases]{parsep=0pt plus 1pt}
\setlist[Cases,1]{wide=0pt, listparindent=\parindent,
    label = \textbf{Case~\arabic*:}, ref = \arabic*}
\setlist[Cases,2]{wide, labelindent=\parindent,
    label = \textbf{Case~\arabic{Casesi}}-(\arabic{Casesii}):}
\setlist[Cases,3]{wide, labelindent=1.5\parindent, topsep=5pt,
    label = \textbf{Case~\arabic{Casesi}}-(\arabic{Casesii})-(\arabic{Casesiii}):}

\crefname{Casesi}{case}{cases}

\NewCommandCopy\statement\equation
\NewCommandCopy\endstatement\endequation

\ExplSyntaxOn

\NewDocumentCommand{\statementapp}{smmmm}
 {
 \IfBooleanF{#1}{\begin{center}}%
  \rainbowrosenfeld_statementapp:nnnn { #2 } { #3 } { #4 } { #5 }%
 \IfBooleanF{#1}{\end{center}}%
 }

\cs_new_protected:Nn \rainbowrosenfeld_statementapp:nnnn
 {
  \begin{tabular}
   {
    |r<{\hspace{0.3pc}}||*{\clist_count:n { #3 }}{>{\hspace{0.1pc}$}c<{$}|}
   }
  \hline
  #1 & \clist_use:nn { #3 } { & } \\
  \hline
  #2 & \clist_use:nn { #4 } { & } \\
  \hline
  \end{tabular}
 }

\ExplSyntaxOff

\NewDocumentCommand{\statementinput}{s m m}{%
    \IfBooleanTF{#1}
        {\statementapp*{Parameters}{Playing the role of}{#2}{#3}}
        {\statementapp{Parameters}{Playing the role of}{#2}{#3}}%
}

\newcommand{\TT}{\mathbf{T}}

\newcommand{\calC}{\mathcal{C}}
\newcommand{\calD}{\mathcal{D}}

\newcommand{\calH}{\mathcal{H}}
\newcommand{\ri}{\mathrm{i}}
\newcommand{\rj}{\mathrm{j}}
\newcommand{\ii}{\iota}

\newcommand{\defeq}{\coloneqq}

\let\originalleft\left
\let\originalright\right
\renewcommand{\left}{\mathopen{}\mathclose\bgroup\originalleft}
\renewcommand{\right}{\aftergroup\egroup\originalright}

\newcommand{\divides}{\mid}

\newcommand{\dirc}{{(\mathrm{d})}} 
\newcommand{\oscl}{{(\mathrm{o})}} 

\newcommand{\ini}{\mathrm{ini}}
\newcommand{\fin}{\mathrm{fin}}

\DeclareMathOperator{\rev}{rev}
\DeclareMathOperator{\shift}{shift}

\newcommand{\Hpart}{\mathbf{H}}

\makeatletter
\DeclareRobustCommand{\vvrev}[1]{\mathpalette\do@vv@rev{#1}}
\newcommand{\do@vv@rev}[2]{%
  \fix@vv@rev{#1}{+}%
  \reflectbox{$\m@th#1\vv{\reflectbox{$\fix@vv@rev{#1}{-}\m@th#1#2\fix@vv@rev{#1}{+}$}}$}%
  \fix@vv@rev{#1}{-}%
}
\newcommand{\fix@vv@rev}[2]{%
  \ifx#1\displaystyle \mkern#23mu
  \else\ifx#1\textstyle \mkern#23mu
    \else\ifx#1\scriptstyle \mkern#22mu
      \else \mkern#22mu
      \fi\fi\fi%
}
\makeatother

\DeclareRobustCommand{\arc}{\vv}
\DeclareRobustCommand{\arcrev}{\vvrev}

\newcommand{\rarc}{\rightarrow}

\newcommand{\outdir}{\Rightarrow}
\newcommand{\indir}{\Leftarrow}

\setcommentmarkup{{\IfIsColored{\color{authorcolor}}{} \textbf{(#3~\arabic{authorcommentcount})} #1}}

\definechangesauthor[name=Jaehyeon, color=ForestGreen]{JS}
\definechangesauthor[name=Jaehoon, color=blue!50!black]{JKim}
\definechangesauthor[name=Deb, color=red]{D}
\definechangesauthor[name=Hyunwoo, color=Mahogany!60!black]{H}

\title{Transversal cycles and paths in tournaments}

\author{
Debsoumya Chakraborti\thanks{
Mathematics Institute, University of Warwick, Coventry, UK. 
Email: {\tt debsoumya.chakraborti@warwick.ac.uk}. 
Supported by the European Research Council (ERC) under the European Union Horizon 2020 research and innovation programme (grant agreement No. 947978). 
}
\and Jaehoon Kim\thanks{
Department of Mathematical Sciences, KAIST, South Korea. E-mail: {\tt jaehoon.kim@kaist.ac.kr}.
Supported by the National Research Foundation of Korea (NRF) grant funded by the Korea government(MSIT) No. RS-2023-00210430.
}
\and Hyunwoo Lee\thanks{Department of Mathematical Sciences, KAIST, South Korea and Extremal Combinatorics and Probability Group (ECOPRO), Institute for Basic Science (IBS). Email:{\tt hyunwoo.lee@kaist.ac.kr}.
Supported by the National Research Foundation of Korea (NRF) grant funded by the Korea government(MSIT) No. RS-2023-00210430 and by the Institute for Basic Science (IBS-R029-C4).}
\and Jaehyeon Seo\thanks{
Department of Mathematics, Yonsei University, South Korea. E-mail: \texttt{jaehyeonseo@yonsei.ac.kr}.
Supported by the National Research Foundation of Korea (NRF) grant funded by the Korea government(MSIT) No.~2022R1C1C1010300 and the Yonsei University Research Fund 2023-22-0125.}}

\date{\today}

\begin{document}
\maketitle

\begin{abstract}
    Thomason~[\emph{Trans.\ Amer.\ Math.\ Soc.} 296.1 (1986)] proved that every sufficiently large tournament contains Hamilton paths and cycles with all possible orientations, except possibly the consistently oriented Hamilton cycle.
    This paper establishes \textit{transversal} generalizations of these classical results. For a collection \(\TT=\{T_1,\dots,T_m\}\) of not-necessarily distinct tournaments on the common vertex set $V$, an $m$-edge directed subgraph $\mathcal{D}$ with the vertices in $V$ is called a transversal if there exists an bijection \(\varphi\colon E(\mathcal{D})\to [m]\) such that \(e\in E(T_{\varphi(e)})\) for all \(e\in E(\mathcal{D})\). We prove that for sufficiently large $n$, there exist transversal Hamilton cycles of all possible orientations possibly except the consistently oriented one. We also obtain a similar result for the transversal Hamilton paths of all possible orientations. 
    These results generalize the classical theorem of Thomason, and our approach provides another proof of this theorem. 
\end{abstract}

\section{Introduction}
This paper studies \textit{transversal} generalizations of the classical results on finding all possible orientations of Hamilton path and cycle in all sufficiently large tournaments. Before diving into the notion of \textit{transversal}, we briefly mention a few of these classical results regarding the existence of paths and cycles in tournaments. 

\subsection{Hamiltonicity in tournaments}
It is well-known that any tournament has a directed Hamilton path. In search of spanning paths with other orientations, Gr\"{u}nbaum~\cite{grunbaum1971antidirected} and Rosenfeld~\cite{rosenfeld1972antidirected} showed that every sufficiently large tournament contains a Hamilton path with an orientation where every two adjacent edges have opposite directions (also, see \cite{hell2002antidirected}). Rosenfeld~\cite{rosenfeld1972antidirected} further conjectured that every sufficiently large tournament contains every orientation of the Hamilton path. Thomason~\cite{thomason1986paths} proved this for tournaments of order at least \(2^{128}\), and Havet and Thomass\'{e}~\cite{havet2000oriented} identified with exactly three exceptions and proved this for all other tournaments.

For the case of Hamilton cycle, a tournament may not contain a directed Hamilton cycle, which is a Hamilton cycle with a consistent orientation. However, it is folklore that any \textit{strongly connected} tournament has a directed Hamilton cycle. Here, a tournament is called strongly connected if, for every pair of vertices $x$ and $y$, there is a directed path from $x$ to $y$. 
For general orientations, in \cite{petrovic1984antidirected,rosenfeld1974antidirected,thomassen1973antidirected}, it was proved that every sufficiently large tournament of even order contains a Hamilton cycle with an orientation where every two adjacent edges have opposite directions. It was also conjectured in \cite{rosenfeld1974antidirected} that every sufficiently large tournament contains every orientation of the Hamilton cycle except possibly the directed one. Thomasson~\cite{thomason1986paths} proved this for tournaments of order at least \(2^{128}\), which was improved to \(68\) by Havet~\cite{havet2000orientedcycles}. 
Recently, Zein claimed a precise result specifying all the exceptional tournaments in \cite{zein2022oriented}.

\begin{theorem}[\cite{thomason1986paths}] \label{thm:Thomason}
    For every sufficiently large $n$, every tournament on $n$ vertices contains all oriented Hamilton paths and cycles, except possibly the directed Hamilton cycle.
\end{theorem}

\subsection{Transversals}
The notion of \textit{transversal} appears in diverse branches of mathematics. A transversal $X$ over a collection $\mathcal{F}=\{F_1,\dots, F_m\}$ of objects stands for an object intersecting every $F_i$. Here, the object can be any mathematical object, such as sets, spaces, set systems, matroids, etc. 
Several classical theorems, which guarantee the existence of certain objects, have been reiterated in the context of transversals, showing the possibility and impossibility of obtaining such an object as a transversal over a certain collection.
For instances, transversals are considered for Carath\'eodory's theorem~\cite{barany1982generalization,kalai2009colorful}, Helly's theorem~\cite{kalai2005topological}, Erd\H{o}s-Ko-Rado theorem~\cite{aharoni2017rainbow}, Rota's basis conjecture~\cite{huang1994relations,pokrovskiy2020rota}, etc.

Following this trend, there has been ample research studying \textit{transversal} generalizations of classical results in graph theory. Other than being interesting in their own right, they often provide a strengthening of the original results. 
Transversals over a collection of graphs have been implicitly used in the literature as the name of rainbow coloring and explicitly defined in~\cite{joos2020rainbow} as follows. 

\begin{definition}
For a given collection $\mathcal{F}=(\calD_1,\dots, \calD_m)$ of graphs/hypergraphs/directed graphs with the same vertex set $V$, an $m$-edge graph/hypergraph/directed graph $\calD$ on the vertex set $V$ is an $\mathcal{F}$-transversal if there exists a bijection $\phi:E(\calD)\rightarrow [m]$ such that $e\in E(\calD_{\phi(e)})$ for all $e\in E(\calD)$.
\end{definition}

The function $\phi$ is often called a \emph{coloring} because one can consider each $\calD_i$ as the set of edges that can be colored with the color $i$.
If the coloring $\phi$ is injective, we say that it is \emph{rainbow} and say the digraph \(\calD\) is a \textit{rainbow} if it is equipped with a rainbow coloring. 
The central question in studying transversal generalizations can be depicted as follows. 

\begin{question}\label{question}
    For a given graph/hypergraph/directed graph $\calD$ with $m$ edges, which properties $\mathcal{P}_n$ will ensure the following? Every collection $\mathcal{F}$ of $m$~graphs/hypergraphs/directed graphs on the same vertex set of size $n$ satisfying the property $\mathcal{P}_n$ contains a $\mathcal{F}$-transversal copy of $\calD$.
\end{question}

Since every object in the family $\mathcal{F}$ can be the same, a necessary condition for a positive answer to the above question is that every $n$-vertex graph/hypergraph/directed graph satisfying the property $\mathcal{P}_n$ contains a copy of~$\calD$. This necessary condition often turns out to be sufficient. 
For instance, the transversal version of Dirac's theorem states that every graph collection $\mathcal{G}=(G_1,\dots,G_n)$ on $n$ vertices satisfying the minimum degree condition $\min_{i\in [n]} \{\delta(G_i)\} \ge \frac{n}{2}$ contains a transversal Hamilton cycle. This was asymptotically established by Cheng, Wang, and Zhao~\cite{cheng2021pancyclicity}, and  completely resolved by Joos and the second author~\cite{joos2020rainbow}. Similar minimum degree conditions for transversal generalizations are proved for other graphs such as $K_r$-factors and bounded degree spanning trees \cite{montgomery2022transversal}, powers of Hamilton cycles \cite{gupta2023general}, and the graphs with bounded degree and sublinear bandwidth~\cite{chakraborti2023bandwidth}. 

However, the transversal versions do not always straightforwardly generalize the classical results. 
For instance, Aharoni, DeVos, de la Maza, Montejano, and \v{S}\'{a}mal~\cite{aharoni2020rainbow} established a transversal version of the classical Mantel's theorem showing that if a graph collection $\mathcal{G}=(G_1,G_2,G_3)$ on a common vertex set of size $n$ satisfies ${\min_{i\in [3]} \{e(G_i)\} > \big(\frac{26-2\sqrt{7}}{81}\big)n^2}$, then it has a triangle transversal, and this condition with the irrational multiplicative constant is best possible. This result shows a different phenomenon compared to previous examples since the number $\frac{26-2\sqrt{7}}{81}$ is larger than $\frac{1}{4}$, which we obtain from Mantel's theorem. 

Transversal generalizations were also considered for hypergraphs \cite{cheng2023rainbow,cheng2021transversal,tang2023rainbow}, directed graphs \cite{babinski2023directed,gerbner2024directed}, tournaments \cite{chakraborti2023hamilton}, and random graphs \cite{anastos2023robust,ferber2022diractype}.

\subsection{Main results}
Previously, in \cite{chakraborti2023hamilton}, we proved the following theorems regarding the existence of a directed Hamilton path and cycle in a family of tournaments as transversal. In this paper, we always assume that \(\TT=\{T_1,\dots,T_m\}\) is a collection of tournaments on the common vertex set $V(\TT)$.

\begin{theorem}[\cite{chakraborti2023hamilton}] \label{thm:rainbow-dir-ham-path-cycle}
For every sufficiently large \(n\), every collection $\TT$ of $n-1$ tournaments with $|V(\TT)|=n$ contains a directed Hamilton path as a transversal.
\end{theorem}
\begin{theorem}[\cite{chakraborti2023hamilton}] 
For every sufficiently large $n$, every collection $\TT$ of $n$ tournaments with $|V(\TT)|=n$ satisfies the following. If all tournaments in $\TT$, possibly except one, are strongly connected, then it contains a directed Hamilton cycle as a transversal.
\end{theorem}

In contrast to the classical non-transversal result, 
it is also shown that the above results fail to hold for small $n$. In this paper, we establish the transversal generalization of \Cref{thm:Thomason},
providing the existence of all transversal cycles and paths in a collection of tournaments. 

\begin{theorem}\label{thm:rainbow-rosenfeld-path}
    For every sufficiently large \(n\), every collection $\TT$ of $n-1$ tournaments with $|V(\TT)| = n$ contains a transversal of every orientation of the Hamilton path.
\end{theorem}

\begin{theorem}\label{thm:rainbow-rosenfeld-cycle}
    For every sufficiently large \(n\), every collection $\TT$ of $n$ tournaments with $|V(\TT)| = n$ contains a transversal of every orientation of the Hamilton cycle except possibly the directed one.
\end{theorem}

Note that our theorems yield almost-pancyclicity in the following strong sense.
For any $m$-vertex set $U$ and a subfamily $\TT'$ of $\TT$ with $m-1$ or $m$ tournaments, if $m$ is sufficiently large, we may apply our results to the smaller collection $\TT'[U]$ of tournaments induced by $U$. 
Then our results show that $\TT$ contains any rainbow cycles and paths with any orientations (possibly except the directed cycle) with any lengths between a large constant and $n$. Moreover, we can even further specify the vertex set of the cycles and paths, and the set of colors used in them.

Our proof strategy is significantly different compared to Thomason~\cite{thomason1986paths} or any other previously used approach, and thus also gives a new proof of \Cref{thm:Thomason}. 
We note that most transversal results for embedding spanning substructures make use of the corresponding non-transversal result in the proof. However, surprisingly, \Cref{thm:Thomason} was not helpful to us in proving \Cref{thm:rainbow-rosenfeld-path,thm:rainbow-rosenfeld-cycle}. Thus, one unique feature of our proof compared to other spanning transversal results in the literature is that we did not use the corresponding non-transversal version in our proof. It will be indeed interesting to find a shorter or conceptually different proof of our main results using \Cref{thm:Thomason} as black box. 

It is worth mentioning that, to prove \Cref{thm:rainbow-rosenfeld-path,thm:rainbow-rosenfeld-cycle}, we utilized a novel way of partitioning the vertices of the tournament, which was introduced in \cite{chakraborti2023hamilton} with the name of \emph{$\Hpart$-partition}.
In this paper, we further develop this concept into a robust version (see \Cref{sec:H-partition}) which helps us to establish an absorption lemma (see \Cref{cor:rainbowDHP_digphs-attached}) for our purpose. We expect this robust $\Hpart$-partition to have further applications for problems dealing with paths in tournaments.

\paragraph{Organization} We organize the rest of the paper as follows. \Cref{sec:prelim} is devoted to notations and some key insights of our proof. In \Cref{sec:tools}, we collect and develop tools that will be useful to us. 
We prove lemmas for vertex absorption in \Cref{sec:absorption} and lemmas for finding near-rainbow brooms in \Cref{sec:embedding_rainbow-brooms}.
We prove \Cref{thm:rainbow-rosenfeld_path_short-blocks} regarding paths with only short blocks in \Cref{sec:shortblocks}, and we prove \Cref{thm:rainbow-rosenfeld-path_long-block-end} regarding paths with a long block in \Cref{sec:rainbow-paths_long-blocks}. These two theorems will directly imply \Cref{thm:rainbow-rosenfeld-path}.
Finally, in \Cref{sec:cycle}, we finish proving \Cref{thm:rainbow-rosenfeld-cycle} using \Cref{thm:rainbow-rosenfeld_path_short-blocks,thm:rainbow-rosenfeld-path_long-block-end}. We end with a few concluding remarks in \Cref{sec:concluding remarks}.

\section{Preliminaries}
\label{sec:prelim}
\subsection{Notations} 

For natural numbers $a\le b$, we write $[a]$ and $[a,b]$ to denote the sets $\{1,\ldots,a\}$ and $\{a,a+1,\ldots,b\}$, respectively. Every \(\log\) in this paper is of base \(2\).
If we say that a claim holds when ${0< \delta \ll \gamma, \beta\ll \alpha<1}$, it means that there exist non-decreasing functions $f : (0,1] \rightarrow (0,1]$ and $g : (0,1]^2 \rightarrow (0,1]$ such that the claim holds when $\gamma,\beta \leq f(\alpha)$ and $\delta \leq g(\gamma, \beta)$. We will not explicitly calculate these functions.
We omit floors and ceilings where they are not crucial.

\paragraph{Digraph} 
Consider digraphs $\calD$ and $\calD'$. We denote the set of vertices and arcs of \(\calD\) by $V(\calD)$ and $E(\calD)$, respectively, and we let $e(\calD) \defeq |E(\calD)|$.
We denote the disjoint union of $\calD$ and $\calD'$ by $\calD \cup \calD'$.
For a vertex $v\in V(\calD)$, $d_\calD^+(v)$ and $d_\calD^-(v)$ are the out-degree and in-degree of $v$, respectively, and $N_\calD^+(v)$ and $N_\calD^-(v)$ are the out-neighborhood and in-neighborhood of $v$, respectively. 
For \(U\subseteq V(\calD)\) and \(\sigma\in\{+,-\}\), we define \(N_\calD^\sigma(v,U)\defeq N_\calD^\sigma(v)\cap U\) and \(d_\calD^\sigma(v,U)\defeq |N_\calD^\sigma(v,U)|\). We often omit the subscript \(\calD\) if it is clear from the context.
For $U\subseteq V(\calD)$, we write $\calD[U]$ for the digraph induced by the vertex set $U$, and write $\calD \setminus U$ or $\calD - U$ for the digraph induced by the vertex set $V(\calD)\setminus U$. When $U=\{u\}$, we sometimes write $\mathcal{D}-u$ instead of $\mathcal{D}-U$.

For vertices \(u\) and \(v\), the arc from \(u\) to \(v\) is denoted by \(\arc{uv}\) or \((uv)^+\) or \(\arcrev{vu}\) or \((vu)^-\). 
Sometimes if the direction of the arc is clear, we omit the arrow sign and write $uv$ to denote the arc. 
We often write \(u\rarc v\) to say that the arc $\arc{uv}$ is present in the digraph of our interest.
For disjoint vertex sets $X$ and $Y$, we denote by $E[X, Y]$ the set of arcs directed from a vertex in $X$ to a vertex in $Y$. 
For a given digraph \(\calD\) and its disjoint vertex subsets \(X\) and \(Y\), we write \(X\outdir Y\) or \(Y\indir X\) when \(\arc{xy}\in E(\calD)\) for every \((x,y)\in X\times Y\).

\paragraph{Oriented path} 
In this paper, a \emph{path} denotes an oriented path with a specified orientation. In other words,
a \emph{path} $P$ is a sequence of vertices $v_1\dots v_{\ell+1}$ with arcs $(v_iv_{i+1})^{\sigma_i}$ with orientations $\sigma_i \in \{+,-\}^{\ell}$ for every $i\in [\ell]$. 

The \emph{length} $\ell(P)$ of \(P\) is the number of arcs in it. 
A path of length $0$ is a single vertex, and we call it an \emph{empty path}. 
If all the orientations $\sigma_i$ of $P$ are the same, then we call it a \emph{directed path}. We write \(\arc{P}_\ell\) and \(\arcrev{P}_\ell\) to denote the directed path of length $\ell$ with the forward and backward direction, respectively. A directed path with forward arcs is often referred to as simply a directed path when there is no ambiguity. 
For a given path $P = v_1 \ldots v_{\ell + 1}$, we denote by \(\rev(P)\) the path \(v_{\ell+1} \dots v_1\) such that \(E(\rev(P))=E(P)\), which is the path we obtain by reversing the order of vertices in \(P\).

A \emph{block} of an oriented path is a maximal directed subpath. We denote by \(b(P)\) the number of blocks of \(P\). 
For two paths $P=v_1\dots v_{\ell+1}$ and $Q=u_1\dots v_{\ell'+1}$ which are vertex-disjoint except $v_{\ell+1}=u_1$,  the path $PQ = v_1\dots v_\ell u_1\dots u_{\ell'+1}$ where \(E(PQ)=E(P)\sqcup E(Q)\) is called a \emph{concatenation} of $P$ and $Q$.
A \emph{decomposition} of the path $P$ is a sequence $P_1,\dots,P_s$ of paths whose concatenation $P_1\dots P_s$ yields $P$. Note that $P_1,\dots, P_s$ are arc-disjoint but two consecutive paths share a vertex. A \emph{block decomposition} of a path $P$ is the decomposition $P_1 \dots P_s$ of $P$ where each $P_i$ is a block of $P$.

\paragraph{Collection of tournaments}
Consider a collection \(\TT\) of tournaments. We write $|\TT|$ to denote the number of tournaments in this collection. Unless stated otherwise, we assume \(\TT\) is of the form \(\{T_c:c\in C\}\) for some set \(C\). We say that this set \(\Gamma(\TT)\defeq C\) is the \emph{color set} of \(\TT\). 
For \(X\subseteq V(\TT)\) and $A\subseteq \Gamma(\TT)$, we define the \emph{vertex-induced collection} $\TT[X] \defeq\{T[X]:T\in\TT\}$ and the \emph{color-induced collection} $\TT_A \defeq \{T_i: i\in A\}$. This naturally yields $\TT_A[X] = \TT[X]_A$. We write $\TT\setminus X$ to denote the induced collection $\TT[V(\TT)\setminus X]$.

Let $\TT = \{T_1, \dots , T_m\}$ be a collection of tournaments. For \(\gamma\in(0,1]\), we define $\TT^{\gamma}$ as the digraph on the vertex set $V(\TT)$ with the arc set
$$  \bigl\{ \arc{uv} :  \left|\{i\in[m]:\arc{uv}\in E(T_i)\}\right|\ge\gamma m \bigr\}.$$
A simple pigeonhole principle implies that \(\TT^{\gamma}\) contains at least one tournament as a subdigraph when $\gamma \leq 1/2$.
We write $\TT_A^{\gamma} = (\TT_A)^{\gamma}$ to make the order of subscripts and superscripts clear.

\paragraph{Coloring}

Consider a color set \(C\),  a digraph \(\calD\) on \(V(\TT)\) and its coloring \(\varphi\) on the arcs of $\calD$. We say \(\varphi\) is \emph{\(C\)-rainbow} if it is injective and its image lies in $C$. If the $C$-rainbow coloring $\varphi$ is clear from the context, we say \(\calD\) is \(C\)-rainbow. 
If in addition the image of \(\varphi\) contains a subset \(C'\subseteq C\), then we say that \(\varphi\) is a \emph{\((C,C')\)-rainbow coloring}.
For a given coloring $\varphi$, we denote \(\varphi(\calD)\) for the set of the colors used by \(\varphi\) to color \(\calD\).  
Unless stated otherwise, a subdigraph \(\calD'\) of \(\calD\) is always considered to have the coloring \(\varphi|_{\calD'}\) inherited from \(\calD\).

\subsection{Proof sketches}
\label{sec:proofsketch}
For finding a transversal copy of a spanning subgraph in $\TT$, we encounter two main difficulties. We have exactly the right number of vertices and exactly the right number of colors.
If there were surplus vertices or colors, then we would have some space to maneuver. However, we do not have such a luxury in our case.

Utilizing the color absorber introduced by Montgomery, M\"{u}yessor, and Pehova \cite{montgomery2022transversal},  the following proof scheme is natural.
Suppose that we want to find a copy of a path $P$ which has a decomposition $P_1 P_2 P_3 P_4$ into four subpaths. \medskip 

\noindent {\bf Step 1. Find a color absorber {\boldmath $(P_1,A,C)$}.}
Consider a tournament  $T\subseteq \TT^{1/2}$ and find a copy of $P_1$ in $T$. As each arc of $T$ has many (at least $\frac{1}{2}|\TT|$) choices of colors, one can choose two disjoint color sets $A,C\subseteq \Gamma(\TT)$ satisfying the following: for any subset $C'\subseteq C$ with $|A\cup C'| = \ell(P_1)$, we can find a $(A\cup C')$-rainbow coloring of $P_1$.
Until the last step, we leave the arcs of $P_1$ uncolored and the colors in $A$ unused. Let $B=\Gamma(\TT)\setminus(A\cup C)$.
\smallskip 

\noindent {\bf Step 2. Find a rainbow copy of {\boldmath $P_2$}.}
We find a rainbow copy of $P_2$ starting from the last vertex of the path $P_1$ in our color absorber. While finding this rainbow path, we only use the colors in $B\cup C$ in such a way that almost all colors in $B$ are used except a very small number of colors.
\smallskip 

\noindent {\bf Step 3. Exhausting all colors in {\boldmath $B$} while embedding a rainbow {\boldmath $P_3$}.}
We find a rainbow copy of $P_3$ starting from the last vertex of our rainbow copy of $P_2$. While finding this rainbow path, we only use the colors in $B\cup C$ in such a way that all colors in $B$ are used. 
\smallskip 

\noindent {\bf Step 4. Vertex absorption.}
We find a rainbow copy of $P_4$ starting from the last vertex of $P_3$ using all vertices using the remaining colors in $C$. As we have to find a spanning copy of $P_4$ in the remaining vertex set, we need to exploit additional tools to deal with this.
\smallskip 

\noindent {\bf Step 5. Color absorption.}
As all colors in $B$ are used and the colors in $A$ are never used, we have the colors in $A\cup C'$ are available with some $C'\subseteq C$ and $\ell(P_1)=|A\cup C'|$. Hence, the property of the color absorber ensures a $(A\cup C')$-rainbow coloring of $P_1$, completing a desired rainbow copy of $P_1$.
\smallskip 

As there are several difficulties following this scheme, we make significant modifications on this scheme.
Below, we explain our ideas for overcoming them.
First, Steps 1 and 5 are standard by now using \Cref{lem:absorber}, which is a reformulation of the color absorbing lemma in \cite{montgomery2022transversal}.

\paragraph{Building paths in Steps 2 and 3.}
Consider Steps 2 and 3. We want to grow a rainbow path into a longer path using colors in $B\cup C$ while forcing the colors in $B$ to be all used, and the final path is the desired $P_i$ for some $i\in \{2,3\}$.
A very natural approach is to grow a path arc by arc, using a new color. However, it is difficult because the orientation of the next arc in the path is arbitrary and may not be compatible with the tournaments of the remaining colors. 
In order to obtain some flexibility, we postpone our decision on the choice of the last arc, by leaving many possibilities. 

Instead of constructing a rainbow path, we find a digraph called \emph{near-rainbow broom}, which is a rainbow path $x_1\dots x_{k-1}$ with a set $Z=\{z_{1},\dots, z_{s}\}$ of vertices where the arc $x_{k-1} z_i$ has the same desired orientation, say $\arc{x_{k-1}z_i}$. In addition, we ensure that the path $x_1\dots x_{k}$ is rainbow no matter which $z_i$ plays the role of $x_{k}$. While the vertices $z_i$ are all candidates for $x_{k}$, we make the decision of choosing $x_k$ only after exploring the candidates for $x_{k+1}, x_{k+2},\dots, x_{k+s'}$ for some $s'<s$. 
To better explain the ideas, assume that $P_i$ is an anti-directed path where every pair of consecutive arcs has different orientations. Choose distinct colors $c,c'$ and assume $\arc{z_1z_2}\in E(T_c)$.
If some $z_i$ has many in-neighbors in $T_{c'}$ for $i\in [2]$, then we take $x_{k-1}z_i$ as the next arc. Otherwise, as $z_1$ has many out-neighbors in $T_{c'}$, then we take $x_{k-1} z_2 z_1$ as the two next arcs. This path, together with the neighbors, yields a new longer broom to continue. 
Moreover, the colors $c$ and $c'$ were arbitrary choices, so we have some additional flexibility in choosing the colors.

We can elaborate on this idea for more general paths. Instead of growing the path arc by arc, we find a path inside the set $Z$ to obtain a new broom with the correct orientation at the end.
This allows us to extend the path as long as the path frequently changes its orientation, i.e. $P_i$ does not have a block much longer than $|Z|$. 
Such arguments can be found in the proofs of \Cref{lem:rainbow-short-free-oscl-double-broom,lem:rainbow-short-lim-oscl-double-broom}.
However, this strategy does not work for the paths $P_i$ with long blocks.
Imagine the situation that the vertices $z_1,\dots, z_s$ have no out-neighbors outside at all in any remaining colors and $P_i$ is a directed path, then there is no way to find a directed path leaving the set $Z$, making $P_i$ impossible to obtain. 
This motivates the definition of an \emph{oscillating path} in \Cref{def:osillating}, which essentially is a path with only very short blocks. 
For oscillating paths, we follow the above idea of brooms.
For directed paths, we utilized the theory we developed before in \cite{chakraborti2023hamilton}. In there, we have introduced a new concept of an $\Hpart$-partition (see \Cref{def:Hpart}) and developed a theory for finding a rainbow directed path. 

These two different approaches deal with oscillating paths and directed paths, but most of the paths are neither. Hence, we decompose the path $P$ into subpaths, all of which are either oscillating or directed, and employ these two different approaches for each piece. This motivates the definition of an \emph{DO-decomopsition} in \Cref{def:DO}.

\paragraph{Vertex absorption in Step 4.}
For Step 4, we again deal with two different cases depending on whether $P_4$ has a long block or not. In this step, the remaining small number of colors in $B$ are all exhausted, and the number of available colors in $C$ is at least a large constant times $\ell(P_4)$.
Thus, one can greedily find a rainbow coloring for any path in $T\subseteq \TT[U]^{1/2}_{C^*}$ for the remaining vertex set $U$ and remaining color set $C^*\subseteq C$. Hence, in the following explanation, we pretend that colors do not matter, and also pretend that finding a desired path $P_4$ in the tournament $T$ is our only goal. 
\smallskip

\noindent \emph{$\Hpart$-absorption}: If $P_4$ contains a long block, then we can utilize the concept of $\Hpart$-partition. However, as the theory developed in \cite{chakraborti2023hamilton} is not strong enough for our purpose, we introduce an additional concept of robustness (see \Cref{def:Hpart}). Using this robustness on $\Hpart$-partition, we prove \Cref{cor:rainbowDHP_digphs-attached} which shows that we can find a desired spanning rainbow directed path while absorbing some vertices outside the $\Hpart$-partition. In fact, \Cref{cor:rainbowDHP_digphs-attached} turns out to be strong enough to do color absorption as well. Thus, Steps 1 and 5 are not necessary in this case, and our proof is different from the proof scheme of five steps introduced before. 
We informally call this \emph{$\Hpart$-absorption}.
In fact, we even use a different decomposition of $P$ which suits our purpose of $\Hpart$-absorption.
\smallskip

\noindent \emph{Transitive absorption}: On the other hand, if $P_4$ does not contain a long block, then it is difficult to categorize what $P_4$ looks like. Hence, we instead focus on the structures in the tournament $T$.
Again, although we do not have apriori knowledge of the structure of $T$, we can utilize the fact that all tournaments contain one universal structure, which is a transitive subtournament of size $\log |V(T)|$. As a transitive tournament has a very specific structure, it is easy to embed any desired path in it. Using this fact and developing some ideas in \cite{thomason1986paths},  we prove \Cref{lem:Ham-path_p-excep_start-end-in-S} which shows that an $m$-vertex transitive subtournament can absorb $O(\log m)$ vertices outside to form any desired path $Q$, as long as $Q$ has many blocks. Since every $n$-vertex tournament contains a $(\log n)$-vertex transitive subtournament, this provides a vertex absorption for $O(\log\log n)$ vertices if $P_4$ has many blocks. Even though this is a very weak absorption, our ideas for Steps~2 and~3 turn out to be strong enough so that we can only leave a large constant number of vertices outside this transitive subtournament (we will fix this transitive subtournament in advance during Step $1$) to absorb in Step~4. We informally call this \emph{transitive absorption}. \medskip

Based on this proof strategy, we consider the two cases depending on the length of the longest blocks in $P$. The following two results deal with these cases, and we will utilize them to obtain \Cref{thm:rainbow-rosenfeld-path,thm:rainbow-rosenfeld-cycle}.

\begin{restatable}{theorem}{ShortBlocks}
\label{thm:rainbow-rosenfeld_path_short-blocks}
  Let $0< 1/n \ll 1$. Let \(\TT\) be a collection of tournaments such that \(|\TT|=|V(\TT)|=n\) and \(\Gamma(\TT)=[n]\) and $P$ be an $n$-vertex path with all blocks having length at most $(\log n)^{1/2}$. 
  Then $\TT$ contains an $[n-1]$-rainbow copy $Q=v_1\dots v_n$ of $P$ in $\TT$ such that $\arc{v_1v_n}\in T_n$.
\end{restatable}

For a given path $P$, denote by $\shift(P)$ the path obtained from $P$ by removing the last arc and appending it at the start of the path in the reverse direction.

\begin{restatable}{theorem}{LongBlock}
\label{thm:rainbow-rosenfeld-path_long-block-end}
    Let \(0<\frac{1}{n}\ll 1\) and \(\TT\) be a collection of tournaments with \(|V(\TT)|=n\) and \(|\TT|=n-1\).
    If \(P\) be a non-directed path of length \(n-1\) whose longest block has length at least \((\log n)^{1/2}\), then $\TT$ contains a rainbow copy of $P$.
    If the last block is a longest block, then we can also find a rainbow copy of $\shift(P)$ having the same start vertex and end vertex as the copy of $P$.
\end{restatable}

Note that \Cref{thm:rainbow-rosenfeld-path} directly follows from these two lemmas.
Additional conditions on the arc $\arc{v_1v_n}$ and on the path $\shift(P)$ will be useful for us to prove \Cref{thm:rainbow-rosenfeld-cycle}.

\section{Tools and auxiliary results} \label{sec:tools}
This section is devoted to collecting a number of helpful results that are motivated by the proof sketch from the last section. Of course, the first one is the following color absorbing lemma, which is a variation of \cite[Lemma~3.3]{montgomery2022transversal}.

\begin{lemma}[{\cite[Lemma~2.5]{chakraborti2023hamilton}}]\label{lem:absorber}
    Let \(0 < 1/n\ll\gamma\ll\beta\ll\alpha \leq 1/2\). Let \(H\) be a digraph with \(\beta n\leq e(H)\leq (\beta+ \frac{1}{2}\gamma ) n\), and \(\TT\) be a collection of tournaments with \(|V(\TT)|=n\) and \(|\TT|=m\ge\alpha n\). Let \(S\) be a copy of \(H\) in \(\TT^{\alpha}\).
        Then there exist disjoint sets \(A,C\subseteq[m]\), with \(|A|=e(H)-\gamma m\) and \(|C| = 10\beta m\) such that the following property holds. Given any subset \(C'\subseteq C\) of size \(\gamma m\), there is a rainbow coloring of \(S\) in \(\TT\) using colors in \(A\cup C'\).
\end{lemma}

\subsection{Median order and structures in tournaments}
In order to build a structures for $\Hpart$-absorption, we want to find a directed path with some additional properties which enables the absorption, i.e., it can be converted to another directed path after adding/deleting some vertices. For this purpose, the concept of median order is useful.

\begin{definition}
    Let \(\calD=(V,E)\) be a digraph. A \emph{median order} of \(\calD\) is a linear order \((v_1,\dots,v_n)\) of \(V\) which maximizes the number \(|\{\arc{v_iv_j}\in E:1\le i<j\le n\}|\) of forward arcs. We write \(\calD=\langle v_1,\dots,v_n\rangle\) if the given ordering is a median order for \(\calD\).
\end{definition}

We mention a few standard properties of median orders in tournaments from~\cite{havet2000median}. For every tournament \(T\) with \(T=\langle v_1,\dots,v_n\rangle\), \(1\le i<j\le n\), and \(I\defeq\{v_i,v_{i+1},\dots,v_j\}\), the following hold.
\begin{proplist}
    \item \(T[I]=\langle v_i,\dots,v_j\rangle\). \label{median_basic_1}
    \item \(|N^+(v_i)\cap I|\ge |N^-(v_i)\cap I|\). In other words, \(d^+(v_i,I)\ge (|I|-1)/2\).\label{median_basic_2}
    \item \(|N^-(v_j)\cap I|\ge |N^+(v_j)\cap I|\). In other words, \(d^-(v_j,I)\ge (|I|-1)/2\).\label{median_basic_3}
\end{proplist}
The two properties \ref{median_basic_2} and \ref{median_basic_3} directly imply the following two propositions.
\begin{proposition}
\label{prop:num-of-small-indeg-vtxs_small}
    For a tournament \(T\) and \(d\ge0\), there are at most \(2d+1\) vertices whose in-degree (resp.\ out-degree) is at most \(d\).
\end{proposition}

\begin{proposition}
\label{prop:median-order-directed-path}
    For a tournament \(T\) with \(T=\langle v_1,\dots,v_n\rangle\), the path $v_1v_2\ldots v_n$ is a directed Hamilton path.
\end{proposition}

The directed path we obtain from the above proposition is good for us as it satisfies the following proposition. 
Roughly speaking, this proposition together with \ref{median_basic_1} asserts that we can freely delete some vertices from the directed path obtained from \Cref{prop:median-order-directed-path} to obtain a shorter directed path, as long as the deleted vertices are pairwise distance at least $4$ away.

\begin{proposition} \label{prop:removal_of_vertex} 
If $T=\langle x_1,x_2,x_3,x_4,x_5\rangle$ is a median order of $T$, then it contains a directed path from $x_1$ to $x_5$ with the vertex set $V(T)\setminus \{x_3\}$.\end{proposition}
\begin{proof}
If $x_2\rightarrow x_4$, then $x_1x_2x_4x_5$ is a desired directed path in $T$. Thus, we can assume $x_4\rightarrow x_2$. Then, by \ref{median_basic_2} and \ref{median_basic_3}, it is easy to check that $x_1\rightarrow x_4$ and $x_2\rightarrow x_5$ in $T$. Consequently, we have $x_1x_4x_2x_5$ as a desired directed path in $T$.
\end{proof}

When we later prove \Cref{thm:rainbow-rosenfeld-cycle}, there is one case that is not covered by either of \Cref{thm:rainbow-rosenfeld_path_short-blocks,thm:rainbow-rosenfeld-path_long-block-end}. This case is when all the arcs in the cycle are consistently oriented except for exactly one arc. The following proposition will be useful in dealing with this exceptional case.

\begin{proposition} \label{prop:modify_hamilton_path}
Let $T$ be a tournament on at least 5 vertices. Then, there are two Hamilton paths $P$ and $P'$ sharing the first vertex such that $P$ is a directed path and every arc in $P'$ except exactly one has the same orientation as $P$.  
\end{proposition}

\begin{proof}
Let $T = \langle v_1,\ldots,v_n \rangle$. By \Cref{prop:median-order-directed-path}, the path $P = v_1\ldots v_n$ is a directed Hamilton path. By \ref{median_basic_3}, there exists $i\in \{n-3,n-2\}$ such that $v_i \rightarrow v_n$. Now, an application of \ref{median_basic_2} implies there exists $j\in [i+1,n]$ such that $v_{i-1} \rightarrow v_j$. Now, in the path $P' = v_1 \ldots v_{i-1} v_j v_{j+1} \ldots v_n v_i v_{i+1} \ldots v_{j-1}$, every arc is forward except $v_n v_i$ which is a backward arc. This finishes the proof.
\end{proof}

The following fact will be useful for us to find brooms in a tournament. Note that Sumner conjectured that $2n-2$ vertices suffice instead of $4n$ in the proposition below, which is shown to be true for sufficiently large $n$ in \cite{kuhn2011sumner}.

\begin{proposition} [\cite{havet2000median}, {Theorem~4}] 
\label{prop:embed_tree}
Every $4n$-vertex tournament contains all $n$-vertex oriented trees. 
\end{proposition}

Also, the following simple proposition is useful. As \Cref{thm:Thomason} is stronger than this, we can instead use \Cref{thm:Thomason}. However, since we aim to provide a different proof of Thomason's result, we rather use the following easy proposition. Since it is quite straightforward, we omit the proof.
\begin{proposition}\label{prop:transitive}
An $n$-vertex transitive tournament contains all Hamilton paths with any orientations.
\end{proposition}

The following proposition is useful to find a desired path inside the set $Z$ described in \Cref{sec:proofsketch} with specified colors.

\begin{proposition} \label{prop:simple-embed-rainbow-path}
For $\ell\in \{1,2\}$, every collection $\TT$ of $\ell$ tournaments on a common vertex set of size five contains a rainbow copy of all oriented paths of length $\ell$.
\end{proposition}
\begin{proof}
If $\ell=1$, this is trivial. 
We may assume that $\ell=2$ and $P=x_1x_2x_3$ satisfies $\arc{x_1x_2}$ and $(x_2x_3)^{\sigma} \in E(P)$ for some $\sigma\in \{+,-\}$.
Take a median order of $T_1=\langle v_1,\dots, v_5\rangle$. 
By \Cref{prop:num-of-small-indeg-vtxs_small}, 
there exists $i\in [2,5]$ satisfying $d^{\sigma}_{T_2}(v_i) \geq 2$. Taking $w\in N^{\sigma}_{T_2}(v_i)\setminus \{v_{i-1}\}$, the path $v_{i-1} v_i w$ yields a desired rainbow copy of $P$. 
\end{proof}

\subsection{Robust \texorpdfstring{\(\Hpart\)}{\textbf{H}}-partition}\label{sec:H-partition}

The following concept of an $\Hpart$-partition was introduced in 
\cite{chakraborti2023hamilton} and proved to be useful for obtaining a rainbow directed path with desired endpoints.
We again make use of this together with the additional `robustness' condition.
\begin{definition}\label{def:Hpart}
Let $0\le \gamma \le 1$ and $r,\ell$ be positive integers and \(T\) be a tournament of order $n$.
A tuple \((W_1,\dots, W_r, w_1,\dots, w_{r-1})\) of disjoint vertex subsets \(W_1,\dots, W_r\subseteq V(T)\) and distinct vertices \(w_1,\dots, w_{r-1}\) in \(V(T) \setminus \bigcup_{i\in [r]}W_{i}\)  is an \textit{$\Hpart(\ell,\gamma)$-partition} if the following hold.
\begin{proplist}
    \item \(\bigcup_{i\in [r]}W_{i}\cup\{w_1,\dots,w_{r-1}\}=V(T)\). \label{H-partition} 
    \item \( \gamma \ell \leq |W_i| \leq \ell \) for each \(i\in [r]\). \label{H-size} 
    \item  \(W_{i}\outdir \{w_{i}\} \outdir W_{i+1}\) for each \(i\in [r-1]\). \label{H-direction} 
\end{proplist}
Moreover, it is \emph{robust} if the following additionally holds.
\begin{enumerate}[label = {\bfseries \Alph{propcounter}\arabic{enumi}}] \setcounter{enumi}{3}
    \item For each $i\in [r-1]$, there exists $i'\leq i$ and $j'>i$ such that $|E[W_{i'},W_{j'}]| \geq \gamma |W_{i'}||W_{j'}| - n$. \label{H-robust}
\end{enumerate}
\end{definition}

Note that a (robust) \(\Hpart(\ell,\gamma_1)\)-partition is a (robust) \(\Hpart(\ell,\gamma_2)\)-partition whenever \(0<\gamma_2\le\gamma_1\le 1\). 
In order to prove the existence of a robust $\Hpart(\ell,\gamma)$-partition, we make use of the following lemma.

\begin{lemma}\label{lem:robust-balanced-vertex}
    Let $T$ be an $n$-vertex tournament with $n \geq 5$. Then, there exists a vertex \(v\) that satisfies $|E[N^-(v), N^+(v)]| \geq \frac{n^2}{25}$ and \(\min\{d^+(v),d^-(v)\}\ge\frac{n}{25}\).
\end{lemma}
\begin{proof}
Count the number of transitive triangles $uvw$ with $u\rightarrow v \rightarrow w$ and $u\rightarrow w$.
We call the vertex $v$ the \emph{center} of the transitive triangle.
As the function $f(x)=\binom{x}{2}$ is convex, there are at least $\sum_{u\in V(T)} \binom{d^+(u)}{2} \geq n \binom{\frac{1}{n}\sum_{u\in V(T)} d^+(u)}{2} = \frac{n(n-1)(n-3)}{8} \geq \frac{n^3}{25}$ transitive triangles. Hence, by the pigeonhole principle, there exists a vertex $v$, which is the center of at least $\frac{n^2}{25}$ transitive triangles. Then it is easy to check that such a vertex $v$ is a desired vertex.
\end{proof}

Now we prove the existence of a robust $\Hpart$-partition.

\begin{lemma} \label{lem:robust-H-partition}
    Let \(0 < \gamma\le \frac{1}{25}\) and \(\ell,n\) be positive integers with \(4\leq \ell\leq n\). Then, every $n$-vertex tournament \(T\) admits a robust \(\Hpart(\ell ,\gamma)\)-partition.
\end{lemma}
\begin{proof}
Note that the statement is clear when $\ell=n$ as \(V(T)\) is a trivial \(\Hpart(\ell ,\gamma)\)-partition, where robustness is vacuous for this partition. Assume that \(\ell \) is the largest possible number such that \(T\) does not have a robust \(\Hpart(\ell ,\gamma)\)-partition. We can assume that $\ell \ge 4$ as otherwise we are done. By the maximality of $\ell$, the tournament $T$ has a robust \(\Hpart(\ell+1 ,\gamma)\)-partition. 
Among all robust \(\Hpart(\ell+1 ,\gamma)\)-partitions, choose one \(P=(W_1,\dots, W_t, w_1,\dots, w_{t-1})\) such that the number of choices of $i\in [t]$ with $|W_i|=\ell+1$ is as small as possible.

If every $i\in [t]$ satisfies \(|W_i| <\ell+1\), then it is also a robust \(\Hpart(\ell ,\gamma)\)-partition, a contradiction. 
Hence, we may assume that there exists $i_0\in [t]$ with $|W_{i_0}|=\ell +1$. 
By \Cref{lem:robust-balanced-vertex}, we can partition \(W_{i_0}\) into three sets \( W^{-}_{i_0}, \{v\}, W^{+}_{i_0}\) such that \( W^{-}_{i_0} \outdir \{v\} \outdir W^{+}_{i_0}\) and
\(|W^-_{i_0}|, |W^+_{i_0}| \geq \frac{1}{25} |W_{i_0}| \geq \gamma (\ell +1)\) and 
\begin{align}\label{eq: robust balanced part}
    |E[W^-_{i_0} , W^{+}_{i_0}]| \geq \frac{1}{25} |W_{i_0}|^2 \geq \gamma |W^{-}_{i_0}||W^{+}_{i_0}|.
\end{align} 
Now consider the partition 
\begin{align*}
    P' &= (W'_1,\dots, W'_{t+1}, w'_1,\dots, w'_{t}) \\
    &= (W_1,\dots, W_{i_0-1}, W^{-}_{i_0},W^{+}_{i_0}, W_{i_0+1},\dots, W_t, w_1,\dots, w_{i_0-1}, v, w_{i_0},\dots, w_{t-1}).
\end{align*}
We claim that this is a robust \(\Hpart(\ell+1 ,\gamma)\)-partition having one smaller number of $i\in [t+1]$ with $|W'_i|=\ell+1$ compared to $P$.
Indeed, since $W_{i_0}$ is partitioned, it is clear that $P'$ contains less number of vertex sets of size $\ell+1$. Also, it is clear that this is a \(\Hpart(\ell+1,\gamma)\)-partition.
Hence, it suffices to check that it is robust.

When \(i=i_0\), by taking \(i'=i_0\) and \(j'=i_0+1\), we can ensure \ref{H-robust} by using \eqref{eq: robust balanced part}. 
When \(i \in [t+1]\setminus\{i_0\}\), we have \(W'_{i} = W_{\tilde{i} }\) for some \(\tilde{i}\in \{i,i-1\}\), then the robustness of \(P\) yields that there are \(i'\leq \tilde{i}\) and \(j'>\tilde{i}\) such that \(|E[W_{i'}, W_{j'}]|\geq \gamma |W_{i'}||W_{j'}| - n\).
If \(\{i',j'\}\) does not contain \(i_0\), then \((i',j')\) also satisfies \ref{H-robust} for \(P'\) and \(i\).
Otherwise, assume \(j'=i_0\) by symmetry. Then we know
\begin{equation} \label{eq:lower bound}
    |E[W_{i'},W_{i_0}]| \geq \gamma |W_{i'}||W_{i_0}| - n \geq \gamma |W_{i'}||W'_{i_0}| + \gamma |W_{i'}||W'_{i_0+1}|  - n.
\end{equation}
On the other hand, we have
\[
    |E[W_{i'},W_{i_0}]|
    = |E[W_{i'},W'_{i_0}]| + |E[W_{i'}, \{w'_{i_0}\}]|+ |E[W_{i'},W'_{i_0+1}]|
    \leq |E[W_{i'},W'_{i_0}]| + |E[W_{i'},W'_{i_0+1}]| + n.
\]
Putting this together with \eqref{eq:lower bound}, we conclude that either 
\[
|E[W_{i'},W'_{i_0}]| \geq \gamma |W_{i'}||W'_{i_0}| -n \;\; \text{or} \;\;
|E[W_{i'},W'_{i_0+1}]| \geq \gamma |W_{i'}||W'_{i_0+1}| -n.
\]
As each of \(W_{i'}, W'_{i_0}\), and \(W'_{i_0+1}\) is a vertex subset in \(P'\), this shows \((i',i_0)\) or \((i',i_0)\) satisfies \ref{H-robust} for \(P'\) and \(i\). Therefore \(P'\) satisfies \ref{H-robust}. 
This contradicts the minimality of \(P\), hence proving the desired lemma.
\end{proof}

This robustness will later help us for $\Hpart$-absorption. However, we will also need to protect the first and last parts of the partition so that the absorption does not modify them. So, we further need to ensure robustness even after removing the first and last partitions. 
To capture that, we introduce the following slight variant of robust partition. 
\begin{definition}\label{def:goodHpartition}
For a tournament \(T\), an \(\Hpart(\ell,\gamma)\)-partition $(U_0,\ldots,U_{r+1},u_0,\ldots,u_r)$ is called a \textit{good \(\Hpart(\ell,\gamma)\)-partition} if $(U_1,\ldots,U_r,u_1,\ldots,u_{r-1})$ is a robust \(\Hpart(\ell ,\gamma)\)-partition of $T\setminus (U_0\cup U_{r+1}\cup \{u_0,u_r\})$.
\end{definition}

Making some further modifications on the robust $\Hpart$-partition, we can obtain a good $\Hpart$-partition as follows.

\begin{corollary}
\label{cor:robust-H-partition-avoid-first-and-last-parts}
    Let \(0 < 1/n \ll \mu, \gamma \le 50^{-7}\). Then any  $n$-vertex tournament $T$ has a good \(\Hpart(\mu n ,\gamma)\)-partition. 
\end{corollary}

\begin{proof}
Apply \Cref{lem:robust-H-partition} on $T$ to find a robust \(\Hpart(\mu n, 1/25)\)-partition $P:=(W_1,\ldots,W_r,w_1,\ldots,w_{r-1})$ satisfying \ref{H-partition}--\ref{H-robust} with $\ell =\mu n$ and $\gamma = 1/25$. Note that $r\ge \frac{n}{\mu n+1} \geq 50^{6}$. Let $p\in [r]\setminus \{1\}$ be the largest index and $q\in [r-1]$ be the smallest index satisfying
\begin{equation}\label{eq:edges with first and last part}
|E[W_1,W_p]| \geq \frac{1}{25} |W_1||W_p| - n \quad \text{and} \quad |E[W_q,W_r]| \geq \frac{1}{25} |W_q||W_r| - n.
\end{equation}
We remark that such \(p\) and \(q\) exist by \ref{H-robust}. Note that if we simply remove the sets $W_1,W_r$ and the vertices $w_1,w_{r-1}$, the remaining partition may no longer satisfy \ref{H-robust}. Instead, we modify the partition $P$ using the following claim. For each $i\in [r-1]\setminus \{1\}$, we define $U_i:=W_i$ and $u_i:=w_i$. 

\begin{claim}\label{claim:splitting W_1 and W_r}
There exist a partition of $W_1$ into $U_0,\{u_0\},U_1$, and of $W_r$ into $U_r,\{u_r\},U_{r+1}$ such that 
\begin{proplist}
    \item $|U_i|\ge \gamma \mu n$ for each $i\in \{0,1,r,r+1\}$, \label{boundary:H-size}
    \item \(U_{i}\outdir \{u_{i}\} \outdir U_{i+1}\) for each \(i\in \{0,r\}\), and \label{boundary:H-direction}
    \item $|E[U_1,U_p]| \geq \gamma |U_1||U_p|$ and $|E[U_q,U_r]| \geq \gamma |U_q||U_r|$. \label{boundary:H-robust}
\end{proplist} 
\end{claim}

Assuming the above claim, we first prove \Cref{cor:robust-H-partition-avoid-first-and-last-parts}.
Since $P$ satisfies \ref{H-partition}--\ref{H-direction}, and~\ref{boundary:H-size}--\ref{boundary:H-direction} hold, both $(U_0,\ldots,U_{r+1},u_0,\ldots,u_r)$ and $P':=(U_1,\ldots,U_r,u_1,\ldots,u_{r-1})$ are \(\Hpart(\mu n ,\gamma)\)-partitions.
To show the robustness of $P'$, we fix $i\in [r-1]$.
Let $n'$ be the number of vertices in $T\setminus (U_0\cup U_{r+1}\cup \{u_0,u_r\})$.

\paragraph{Case 1. $p \le i < q$.}
Since $P$ satisfies \ref{H-robust}, by the choices of $p,q$, there exists $2\le i'\le i$ and $i< j' \le r-1$ such that $|E[U_{i'},U_{j'}]| = |E[W_{i'},W_{j'}]| \geq \frac{1}{25} |W_{i'}||W_{j'}| - n \geq \gamma |U_{i'}||U_{j'}| - n'$.

\paragraph{Case 2. $i < p < r$.}
By the first inequality of \ref{boundary:H-robust}, with $i':=1\le i$ and $j':=p > i$, we have $|E[U_{i'},U_{j'}]| \geq \gamma |U_{i'}||U_{j'}| - n'$.

\paragraph{Case 3. $i \ge q$.} 
We first note that this case contains the case $p=r$ because $p=r$ implies $q=1$. By the second inequality of \ref{boundary:H-robust}, with $i':=q\le i$ and $j':=r > i$, we have $|E[U_{i'},U_{j'}]| \geq \gamma |U_{i'}||U_{j'}| - n'$. This proves the robustness of $P'$ and
finishes the proof of the Corollary.
\end{proof}
The following proof of \Cref{claim:splitting W_1 and W_r} completes the proof of \Cref{cor:robust-H-partition-avoid-first-and-last-parts}.
\begin{proof}[Proof of \Cref{claim:splitting W_1 and W_r}]
By the first inequality in \eqref{eq:edges with first and last part} and the assumption $1/n\ll \mu,\gamma$, there is a vertex set $V_1\subseteq W_1$ of size $|W_1|/50$ such that every $v\in V_1$ satisfies $d^+_T(v,W_p) \geq |W_p|/50$. Apply \Cref{lem:robust-balanced-vertex} on $T[V_1]$ to find a vertex $u_0$ such that $d^{\sigma}_{T(V_1)}(u_0)\ge |V_1|/25$ for every $\sigma\in \{+,-\}$. Now, we set $U_0:=N^-_{T(W_1)}(u_0)$ and $U_1:=N^+_{T(W_1)}(u_0)$. 
Thus, for every $i \in \{0,1\}$, \ref{boundary:H-size} holds since 
\[
|U_i|\ge \frac{|V_1|}{25} = \frac{|W_1|}{25\cdot 50} \ge \frac{\mu n}{25^2\cdot 50} \ge \gamma \mu n.
\]
By construction, \ref{boundary:H-direction} holds for $i=0$. Finally, the first inequality in \ref{boundary:H-robust} holds since $\gamma < 50^{-7}$ and
\[
|E[U_1,W_p]|= |E[N^+_{T(V_1)}(u_0),W_p]| \ge \frac{|V_1|}{25} \cdot \frac{|W_p|}{50} = \frac{|W_1|}{25\cdot 50} \cdot \frac{|W_p|}{50} \ge \frac{|U_1||W_p|}{25\cdot 50^2}.
\] 
Considering two cases of $(p,q)=(r,1)$ and $(p,q)\neq (r,1)$, 
we see that in either case, the above inequality together with the second inequality in \eqref{eq:edges with first and last part} implies \[
|E[U_q,W_r]| \geq \frac{1}{25\cdot 50^2} |U_q||W_r|.
\]
Next, by the above inequality, there is a set $V_r\subseteq W_r$ containing $|W_r|/50^3$ vertices such that every $v\in V_r$ satisfies $d^-_T(v,U_q) \geq |U_q|/50^3$. 
Apply \Cref{lem:robust-balanced-vertex} on $T[V_r]$ to find a vertex $u_r$ such that $d^{\sigma}_{T(V_r)}(u_r)\ge |V_r|/25$ for every $\sigma\in \{+,-\}$. Now, we set $U_r:=N^-_{T(W_r)}(u_r)$ and $U_{r+1}:=N^+_{T(W_r)}(u_r)$. 
Thus, for every $i \in \{r,r+1\}$, \ref{boundary:H-size} holds since 
\[
|U_i|\ge \frac{|V_r|}{25} = \frac{|W_r|}{25\cdot 50^3} \ge \frac{\mu n}{25^2\cdot 50^3} \ge \gamma \mu n.
\]
By construction, \ref{boundary:H-direction} holds for $i=r$. Finally, again considering two cases of $(p,q)=(r,1)$ or not, the second inequality in \ref{boundary:H-robust} holds since 
\[
|E[U_q,U_r]|= |E[U_q,N^-_{T(V_r)}(u_r)]| \ge \frac{|V_r|}{25} \cdot \frac{|U_q|}{50^3} = \frac{|W_r|}{25\cdot 50^3} \cdot \frac{|U_q|}{50^3} \ge \frac{|U_q||U_r|}{25\cdot 50^6} \ge  \gamma |U_q||U_r|.
\]
This proves the claim.
\end{proof}

In \cite{chakraborti2023hamilton}, we proved the following lemma showing that an $\Hpart$-partition is useful for finding a directed Hamilton path.

\begin{lemma}[{\cite[Lemma~3.4]{chakraborti2023hamilton}}]
\label{lem:rainbowDHP_1st}
    Let \(0<\frac{1}{n}\ll\mu\ll \gamma, \alpha\le 1\). Let \(\TT\) be a collection of tournaments with \(|V(\TT)| = n\) and \(|\TT| = n-1\). 
    Let \(T\subseteq \TT^{\alpha}\) be a tournament. Let \(w_0,w_r\in W\), and \((W_1,\dots,W_r,w_1,\dots,w_{r-1})\) be an \(\Hpart(\mu n,\gamma)\)-partition of \(T\setminus\{w_0,w_r\}\).
    Suppose \(w_0\outdir W_1\) and \(W_r\outdir w_r\) in~$T$.
    Then, there is a rainbow directed Hamilton path $P$ in \(\TT\) from \(w_0\) to \(w_r\).
\end{lemma}

\section{Vertex absorption}
\label{sec:absorption}
In this section, we establish our key lemmas to execute the vertex absorption step as described in the proof sketch of \Cref{sec:proofsketch}. As mentioned there, we will use two different types of vertex absorption.
\subsection{Transitive absorption}
In \cite{thomason1986paths}, Thomason proved the following lemma asserting that a transitive subtournament $S$ can `absorb' a certain small vertex set $W$ outside, by finding a directed path containing $W$ starting from $S$ and ending at $S$.

\begin{lemma}[{\cite[Lemma~8]{thomason1986paths}}]
\label{lem:dir-path_start-end-in-S}
    Let \(S\) be a transitive subtournament of a tournament \(T\) with \(|S|\ge 4\) and $W$ be the vertex set $V(T)\setminus V(S)$. Suppose that every vertex in $W$ has at least one in-neighbor and out-neighbor in $V(S)$. Then for any integer $m$ with \(|W|+3\le m\le |V(T)|-1\), there is an arc \(\arc{uv}\) in \(S\) and a directed path \(P\) of length \(m\) from \(u\) to \(v\) with $W\subseteq V(P)$ such that \(u\outdir (V(T)\setminus V(P))\outdir v\).
\end{lemma}

Developing the idea of Thomason used in \cite[Lemma~10]{thomason1986paths}, we can prove the following lemma showing that a large transitive subtournament is an appropriate absorber for finding a spanning path with sufficiently many blocks.

\begin{lemma}
\label{lem:Ham-path_p-excep_start-end-in-S}
For an integer $p\geq 0$, let $P$ be an $n$-vertex path with $n>2^{p+4}$ having at least $5p$ blocks and $T$ be an $n$-vertex tournament containing a transitive subtournament $S$ on $n-p$ vertices.
Then $T$ contains a copy of $P$ whose both endpoints lie in $S$.
\end{lemma}
\begin{proof}
If $p=0$, then this follows from \Cref{prop:transitive}. Assume $p\geq 1$ and that the statement holds for $p-1$.

First, suppose that \(P\) has a block of length \(\ell \ge p+3\). Without loss of generality, assume that this block is forwardly directed. Let \(W\subseteq V(T)\setminus V(S)\) be the set of vertices having at least one in-neighbor and out-neighbor in $S$.
If \(W=V(T)\setminus V(S)\), then we apply \Cref{lem:dir-path_start-end-in-S} to find an arc \(\arc{uv}\) in \(S\) and a directed path \(Q\) of length \(\ell \) from $u$ to $v$ such that \(V(T)\setminus V(Q)\subseteq S\) and \(u\outdir (V(T)\setminus V(Q))\outdir v\). 
Since \(T\setminus Q\) is a transitive tournament, \Cref{prop:transitive} implies that it contains \(P-V(Q)\), 
which is either a path or two vertex-disjoint paths. Since \(u\outdir (V(T)\setminus V(Q))\outdir v\), this yields a copy of $P$ with both endpoints inside $V(T)\setminus V(Q) \subseteq S$.

If \(V(T\setminus S)\setminus W\) is nonempty, we choose \(w\in V(T\setminus S)\setminus W\). Without loss of generality, assume \(w\indir S\). 
Consider the block decomposition $P_1P_2\dots P_k$ of $P$, then $k\geq 5p$.
Let $x$ be the first (left-most) vertex of in-degree two in $P$ and $y$ the last vertex of in-degree two in $P$. By the definition of the block decomposition, $x$ is the last vertex of $P_{i_*}$ for some $i_*\in \{1,2\}$ and $y$ is the first vertex of $P_{j_*}$ for some $j_*\in \{k-1,k\}$. As $k\geq 5$, we know $i_*< j_*$ and at least one, say $Q_1$, of the paths 
$Q_1= P_{i_*+1}\dots P_{k} - x$ and $Q_2= P_1\dots P_{j_*-1} - y$ has length at least $\frac{1}{2}\ell(P) \ge 2^{p+3}$.
Also, both \(Q_1\) and \(Q_2\) have at least $5p-3\geq 5(p-1)$ blocks.

Since $|V(Q_1)|> 2^{p+3}\geq p$ and \(|V(S)|\ge n-p\), we have \(|V(S)|> \ell(P_1)+\dots + \ell(P_{i_*})\). We delete $\ell(P_1)+\dots + \ell(P_{i_*})-1$ arbitrary vertices in $S$ and the vertex $w$ from $T$ to obtain $T'$. Note that it is easy to see that the remaining part of $S$ still induces a nonempty transitive subtournament.
We apply the induction hypothesis to find a copy $R$ of $Q_1$ in $T$ which starts and ends in $S$, containing all vertices in $V(T')\setminus S = V(T)\setminus (S\cup \{w\})$.
As $S\setminus R$ is a transitive tournament, using \Cref{prop:transitive}, we can find a copy $R'$ of $P_1\dots P_{i_*}-x$ in $S\setminus R$. Because $V(R')\outdir w \indir V(R)$, the concatenation $R' \arc{uw}\arcrev{wv} R$ yields a copy of $P$ where $u$ is the last vertex of $R'$ and $v$ is the first vertex of $R$.

Second, suppose that every block of \(P\) has length at most \(p+2\). 
We decompose $P$ into paths \(Q_1\) and \(Q_2\) and an arc $e$ such that $P=Q_1eQ_2$ and $\ell(Q_1) = 4p+11$. 
Without loss of generality, assume that $e=\arc{xy}$ is forwardly directed with $x\in Q_1$ and $y\in Q_2$.
We partition $S$ into two transitive tournaments $S_1$ and $S_2$ with $S_1 \outdir S_2$ and $|S_1|=4p+11$. (If the arc $e$ was $\arcrev{xy}$, then we swap the roles of $S_1$ and $S_2$ by choosing $|S_2|=4p+11$.)

We choose a vertex $w\in V(T)\setminus V(S)$ and consider the median order $\langle w_1,\dots w_s \rangle$ of $S_1$. 
Consider the tournaments 
\[
    T_1 \defeq T[V(S_1)\cup \{w\}]
    \quad\text{and}\quad
    T_2 \defeq T-(V(S_1)\cup \{w\}).
\]
As \(\ell (Q_2)\ge 2^{p+4}-(4p+12)\ge 2^{p+3}\) and \(b(Q_2)\ge \frac{2^{p+3}}{p+2} > 5(p-1)\), the induction hypothesis yields a copy of $Q_2$ in $T_2$ with both endpoints in $S_2$.
As $S_1\outdir S_2$, it suffices to find a copy of $Q_1$ in $T_1$ with both endpoints in $S_1$. 

Note that at least two arcs from $w$ to $\{w_{2p+5}, w_{2p+6}, w_{2p+7}\}$ has the same orientation. Without loss of generality, assume that $\arc{w w_{i'}}, \arc{w w_{j'}} \in T$ for $i',j'\in [2p+5,2p+7]$.
As $Q_1$ has $4p+11$ arcs and every block has length at most $p+2$, the block decomposition $R_1\dots R_k$ of $Q_1$ satisfies $k\geq 5$. 
As $k\geq 5$, we can find $s\in [k-1]$ such that the last vertex in $R_s$ has in-degree two in $Q_1$. 
As $\ell(R_s)+\ell(R_{s+1})\leq 2p+4 < \min\{i',j'\}$, we can find two vertex-disjoint directed paths one from $w_1$ to $w_{i'}$ of length $\ell(R_s)-1$ and another one from $w_2$ to $w_{j'}$ of length $\ell(R_{s+1})-1$ in the transitive tournament $\langle w_1,\dots, w_{2p+7}\rangle$.
Deleting the vertices in those paths from $S_1$ leaves a transitive tournament, and \Cref{prop:transitive} implies that it contains $Q_1 - R_s-R_{s+1}$. As every vertex in $S_1$ is an out-neighbor of both  $w_1$ and $w_2$, the copy of $Q_1-R_s-R_{s+1}$ together with those paths of length $\ell(R_s)-1, \ell(R_{s+1})-1$ and the vertex $w$ yields a desired copy of $Q_1$. Together with $Q_2$, this yields a desired copy of $P$ in $T$.
\end{proof}


\subsection{\texorpdfstring{\(\Hpart\)}{\textbf{H}}-absorption}

Now we prove the following lemma showing that a robust $\Hpart$-partition of $W$ can absorb a small set of vertices if those vertices have many out-neighbors and in-neighbors in $W$ for enough number of colors.

\begin{lemma}
\label{cor:rainbowDHP_digphs-attached}
    Let \(0<\frac{1}{n}\ll\mu\ll\alpha, \gamma, \varepsilon, \frac{1}{p+1}\le 1\) and
    \(\TT\) be a collection of tournaments. Let \(U, W\subseteq V(\TT)\) be disjoint vertex sets of size \(p\) and \(n\), respectively and
    \(C\subseteq\Gamma(\TT)\) be a color set of size \(n+p-1\). 
    For each $u\in U$, let $D_u\subseteq C$ be a set of at least $3p$ colors satisfying
    \begin{equation} \label{minimum_degree_condition}
        d_{T_d}^+(u,W),\, d_{T_d}^-(u,W)\ge\varepsilon n
        \quad\text{for each }d\in D_u.
    \end{equation}
    Let \(T\subseteq \TT^{\alpha}_C[W]\) be a tournament, $w_0,w_{r}\in W$ be two vertices, and
    \(\calH=(W_1,\dots,W_{r},w_1,\dots,w_{r-1})\) be a robust \(\Hpart(\mu n,\gamma)\)-partition of \(T\setminus\{w_0,w_{r}\}\) such that $w_0\outdir W_1$ and $W_r\outdir w_r$ in $T$.  Then there is a $C$-rainbow directed Hamilton path in \(\TT[U\cup W]\) from $w_0$ to $w_r$.
\end{lemma}

\begin{proof} 
We will take a directed path $P$ in $T[W]$.
Using three different ways, we will absorb the vertices in $U$ into $P$ and choose colors for some arcs to obtain precolored subpaths. We then use \Cref{lem:rainbowDHP_1st} to find paths connecting the precolored subpaths, obtaining a desired rainbow path.
Depending on how we absorb the vertices, we plan to categorize the vertices in $U$.

\paragraph{Step 1. Partitioning $U$ based on the arcs between them and $W$.}
In the first step, we will partition the vertex set $U$ into three sets \(U_1\), \(U_2\), and \(U_3\). Those vertices in different sets will be absorbed in different ways in Step~2.
We plan to take a spanning path $P_i$ in $T[W_i]$ for each $i\in [r]$, and use these paths to absorb some vertices in $U$. In order for this, we need to specify an index $i$ for each vertex $u\in U$ such that the path $P_i$ can absorb the vertex $u$. For this,
consider a median order $\langle v_1^i,v_2^i,\ldots,v_{n_i}^i\rangle$ of the vertices in $T[W_i]$ for each $i\in [r]$, where $n_i\defeq|W_i|$. 
Let $P_i$ be the directed path $v_1^i v_2^i  \cdots v_{n_i}^i$ obtained from the median order. 
    
Let $P$ be the directed path $w_0 P_1 w_1\dots w_{r-1} P_{r} w_r$ in $T$ obtained by concatenating the paths $P_1,\dots, P_{r}$ together with the vertices $w_0, w_1, \dots, w_{r-1}, w_r$. In particular, $P$ is a directed Hamilton path of $T[W]$.
Denote the \(i\)th vertex of \(P\) by \(v_i\) so that $P = v_1\cdots v_{n}$.
    
For each $u\in U$, we choose a set of three colors $D'_u\subseteq D_u$ so that $D'_u\cap D'_{u'}=\emptyset$ for all $u\neq u'\in U$. Indeed this is possible as $|D_u|\geq 3p=3|U|$ for all $u\in U$.

We claim that for each $u\in U$, we can choose $ \ii(u)\in [n-1]$ so that it satisfies one of the following two conditions. 

\begin{proplist}[series=U]
    \item There exist distinct $d_1=d_1(u)$ and $ d_2=d_2(u)$ in $D'_u$ such that $v_{\ii(u)}\in N^-_{T_{d_1}}(u)$ and $v_{\ii(u) + 1}\in N^+_{T_{d_2}}(u)$. \label{U1}
    \item There exist distinct $d_1=d_1(u)$ and $ d_2=d_2(u)$ in $D'_u$ such that $N^+_{T_{d_1}}(u)\supseteq \{v_1,\dots, v_{\ii(u)}\}$ and $N^-_{T_{d_2}}(u)\supseteq \{v_{\ii(u)+1},\dots, v_{n-1}\}$. \label{U23}
\end{proplist}
Indeed, suppose for some $u\in U$, \ref{U1} does not hold.
Then, consider the smallest $i\in [n]$ such that there exists $d\in D'_u$ satisfying $v_i\in N^-_{T_d}(u)$. 
Let $ \{d'_1,d'_2\}\defeq D'_u\setminus \{d\}$. Since \ref{U1} does not hold with the choice of $(d_1,d_2)=(d,d'_j)$, we have $v_{i+1}\in N^-_{T_{d'_j}}(u)$ for each $j\in [2]$. Applying the same reasoning repeatedly, we have $v_{k}\in N^-_{T_{d'_j}}(u)$ for each $j\in [2]$ and $i+1\le k\le n-1$. By the choice of $i$, we have $v_{k}\in N^+_{T_{d'_j}}(u)$ for each $j\in [2]$ and $1\le k\le i-1$. Thus, depending on whether $v_i\in N_{T_{d_j'}}(u)$ or not for each $j\in \{1,2\}$, \ref{U23} holds with a choice of $\ii(u)\in \{i-1,i\}$ and $\{d_1,d_2\}=\{d'_1,d'_2\}$, which proves our claim. 

Furthermore, by \eqref{minimum_degree_condition} and the fact that $\mu \ll \varepsilon$, we have the following for every $u\in U$.
\begin{enumerate}[resume*=U]
    \item If $u\in U$ satisfies \ref{U23}, then
    $n_1< \ii(u)< n - n_r$.  \label{separation from sides}
\end{enumerate}

Let $U_1$ be the set of vertices $u\in U$ where $\ii(u)$ satisfies \ref{U1}, and for each $u\in U_1$, we fix $d_1(u), d_2(u)$ satisfying \ref{U1}.
For the vertices $u \notin U_1$ where $\ii(u)$ satisfies \ref{U23}, we will divide the vertices into two sets $U_2$ and $U_3$.
To define $U_2$ and $U_3$, let 
\begin{equation} \label{def of M}
    M \defeq \bigcup_{i\in [r]} \{ v^i_{ \gamma n_i/100} ,\dots, v^i_{(1-\gamma /100)n_i} \}
\end{equation}
be the set of vertices in a middle part of the paths $P_1,\dots, P_{r}$.
Let $U_2$ be the set of vertices $u\in U$ satisfying \ref{U23} with $v_{\ii(u)} \in W\setminus M$, and $U_3$ be the set of vertices $u\in U$ satisfying \ref{U23} with $v_{\ii(u)}\in M$.

\paragraph{Step 2. Absorbing the vertices in $U$.} 
In this step, we will define index sets \(J_1\), \(J_2\), and \(J_3\), and a directed path $R_j$ for each $j\in J\defeq J_1\cup J_2\cup J_3$. We plan to replace the arc $\arc{v_jv_{j+1}}$ by the directed path $v_j R_j v_{j+1}$ for each $j\in J$. By carefully constructing the paths $R_j$, we will make sure $\bigcup_{j\in J} V(R_j)$ contains all vertices in $U$. The paths $R_j$ will also use some vertices from other parts of $P$ leaving some damages on other parts of the path. However those damages will be well-spread so that we can use \Cref{prop:removal_of_vertex} to recover the path. 
We will use different types of paths $R_j$ for each of $J_1$, $J_2$, and $J_3$.

\subparagraph{Step 2-1. Absorbing \(U_1\).}
In order to absorb vertices in $U_1$, we define 
\[
    J_1 \defeq \{ \ii(u): u\in U_1\}
    \quad\text{and}\quad
    I_1 \defeq \{ i\in [r] : \exists j\in J_1 \text{ such that } v_j \in \{w_{i-1}\}\cup W_i \}.
\]
For each $j\in J_1$, consider the set $U(j)$ of all vertices $u\in U_1$ with $\ii(u)= j$.  Then, the subtournament $T[U(j)]$ contains a directed Hamilton path $R_j$. Let the first vertex and the last vertex of the directed path $R_j$ be $x_j$ and $y_j$, respectively.
For each $j\in J_1$, replace the arc $\arc{v_{j}v_{j+1}}$ in $P$ with the directed path $v_{j} R_j v_{j+1}$. Using \ref{U1}, we color the arc $\arc{v_j x_j}$ with the color $d_1(x_j)$ and the arc $\arc{y_j v_{j+1}}$ with the color $d_2(y_j)$.

Let $P^1$ be the resulting partially colored path.
Note that $|I_1|\leq p$ and the path $P^1$ contains $P_i$ as a subpath for all $i\notin I_1$.
Also, all the uncolored arcs of $P^1$ are in $T$.

The above replacements do not damage the path as we simply insert a path replacing an arc. However, absorbing the vertices in $U_2$ and $U_3$ will remove some vertices from the path $P$.
As we wish to use \Cref{prop:removal_of_vertex}, we need to make sure that the removed vertices are far apart in the path $P$.

We say a set $A$ of natural numbers is \emph{$r$-separated} if $|a-b|>r$ for all distinct $a,b\in A$.
We say a set $X\subseteq V(P)$ is \emph{$r$-separated} if $|a-b|>r$ for all distinct $v_a, v_b\in X$.
For a given vertex set $X\subseteq V(P)$, let $B^r(X) \defeq \{ v_b\in V(P): |b-a|\leq r \text{ for some } v_a\in X\}$.

\subparagraph{Step 2-2. Absorbing \(U_2\).}
In order to absorb a vertex $u$ in $U_2$, we plan to find an arc $v_{\rj(u)}v_{\rj'(u)}$ with $\rj(u)< \ii(u) < \rj'(u)$ and let $R_{\rj(u)}$ be the directed path $\arc{v_{\rj'(u)}u}$ so we can
replace the arc $\arc{v_{\rj(u)}v_{\rj(u)+1}}$ with $v_{\rj(u)} v_{\rj'(u)} u v_{\rj(u)+1}$. 
We want to find those indices using the property \ref{H-robust} in such a way that the vertex set $\bigcup_{u \in U_2} \{v_{\rj(u)}, v_{\rj'(u)}\}$ becomes a $4$-separated set. The following claim ensures such choices.

\begin{claim}\label{claim:choices for U2}
There exist choices $\rj(u),\,\rj'(u)\in [n-1]$ for every $u\in U_2$ such that the following hold.
\begin{enumerate}
    \item $\{ \rj(u),\, \rj'(u): u\in U_2\}$ are all distinct numbers such that $\rj(u)< \ii(u) < \rj'(u)$ for every $u\in U_2$.
    \item For each $u\in U_2$, the arc $\arc{v_{\rj(u)}v_{\rj'(u)}}$ is in $T$.
    \item $J_2\cup J'_2$ is $4$-separated and disjoint from $B^4(J_1)$, where
    \[
        J_2\defeq \{\rj(u) :u\in U_2\}
        \quad\text{and}\quad
        J'_2\defeq \{ \rj'(u):u\in U_2\}.
    \]
    \item For each $u\in U_2$, the vertices $v_{\rj(u)}$ and $v_{\rj'(u)}$ are in $M$ (as defined in \eqref{def of M}).
\end{enumerate}
\end{claim}

\begin{claimproof}
For each $u\in U_2$, without loss of generality, suppose $v_{\ii(u)} \in \{w_{i-1}, v^i_1,\dots, v^i_{\gamma n_i/100}\}$ for some $i\in [r]$.
By \ref{separation from sides}, we must have $i\in [2,r-1]$. As $(W_1,\dots, W_r, w_1,\dots, w_{r-1})$ is a robust $\Hpart(\mu n, \gamma)$-partition, the robustness implies that there exist \(i_1, i_2\in [r]\) such that $i_1\leq i-1 < i_2$ and
$|E[W_{i_1},W_{i_2}]| \geq \gamma n_{i_1}n_{i_2} - n$.
Then, for $A\defeq W_{i_1}\cap M$ and $B\defeq W_{i_2}\cap M$, we have
$$|E(A,B)| \geq |E[W_{i_1},W_{i_2}]| - \frac{4\gamma}{100} n_{i_1}n_{i_2} \geq  \mu^3 n^2.$$
So, there are at least $\mu^3 n^2$ choices of arcs $\arc{v_{\rj(u)}v_{\rj'(u)}}$ for each $u\in U_2$ with $\rj(u)< \ii(u) < \rj'(u)$ and $v_{\rj(u)}, v_{\rj'(u)}\in M$. 
The inequality $\ii(u) < \rj'(u)$ holds because either $i_2>i$ or $v_{\rj'(u)}$ lies in $M$ while $v_{\ii(u)}\in \{ w_{i-1},v^i_1,\dots, v^i_{\gamma n_i/100}\}$.
Note that if $v_{\ii(u)} \in \{ v^i_{(1-\gamma^2)n_i},\dots, v^i_{n_i}\}$ for some $i\in [r]$, then by a similar argument, we can also prove that there are at least $\mu^3 n^2$ choices of such arcs.

For every $i\in J_1$, there are at most $9 \cdot n$ other arcs whose at least one endpoint is in $B^4(v_{\ii(u)})$. Similarly, once we choose an arc $\arc{v_{\rj(u)}v_{\rj'(u)}}$ for a vertex $u\in U_2$, there are at most $2\cdot 9 \cdot n = 18 n$ other arcs whose at least one endpoint is in $B^4(v_{\rj(u)})\cup B^4(v_{\rj'(u)})$. Hence, as $\mu^3 n^2 > 27 n p \ge 9 n |J_1| + 18 n |U_2|$, we can greedily choose an arc $\arc{v_{\rj'(u)}v_{\ri(u)}}$ for each $u\in U_2$ one by one so that $J_2\cup J'_2$ is $4$-separated and disjoint from $B^4(J_1)$. This proves the claim.
\end{claimproof}
    
We take the sets $J_2$ and $J'_2$ as in \Cref{claim:choices for U2}. 
For each $j'\in J'_2$, we apply \Cref{prop:removal_of_vertex} to replace the part $v_{j'-2}v_{j'-1}v_{j'} v_{j'+1}v_{j'+2}$ of the path $P^1$ to a shorter path $v_{j'-2}yy'v_{j'+2}$ in $T$ with $\{y,y'\}=\{v_{j'-1}, v_{j'+1}\}$. Let $P'$ be the resulting path, then $V(P') = V(P^1)\setminus\{ v_{j'}: j'\in J'_2\}$. 
For each $j\in J_2$, there exists unique $u\in U_2$ with $\rj(u)=j$. Then we replace the arc $\arc{v_{j} v_{j+1}}$ with the directed path $v_{j} v_{\rj'(u)} u v_{j+1}$.
Using \ref{U23} and the fact $j(u)<\ii(u)<j'(u)$, we color the arc $\arc{v_{\rj'(u)} u}$ with the color $d_2(u)$ and the arc $\arc{u v_{j+1}}$ with the color $d_1(u)$.

By doing this for all $j\in J_2$, we obtain a path $P^2$ with $V(P^2) = V(P^1)\cup U_2$. Note that $P^2$ is a partially colored path where all the uncolored arcs of $P^2$ are in $T$.
Moreover, by our construction, for every $j\in [n']\setminus B^4(J_1\cup J_2\cup J'_2)$, the path $v_{j-2}v_{j-1}v_jv_{j+1}v_{j+2}$ is still a subpath of $P^2$.
Again, we define 
\[
    I_2 \defeq \{i\in [r]: \exists j \in J_2 \text{ such that } v_j \in W_i \}
    \quad\text{and}\quad
    I'_2 \defeq \{i \in [r]: \exists j \in J'_2 \text{ such that } v_j \in W_i  \}.
\]

\subparagraph{Step 2-3. Absorbing \(U_3\).}
We next absorb the vertices in $U_3$. Let
\[
    I_3 \defeq \{ i\in [r]: v_{\ii(u)} \in W_i \text{ for some }u\in U_3\}.
\]
For each $i\in I_3$, let $Y_i$ be the set of all vertices $u\in U_3$ where $\ii(u)\in M\cap W_i$.
For each $i\in I_3$, we will find a subpath $R_i$ in $W_i$ which starts at a vertex $v_j$ and ends at a vertex $v_{j'}$ for some $j$ and $j'$ where all vertices in $M\cap W_i$ are after $v_j$ and before $v_{j'}$ in the median order of $P_i$. Let $R'_i$ be a directed Hamilton path in $T[Y_i]$, then we can replace the arc $v_{j}v_{j+1}$ with the directed path $v_j R_i R'_i v_{j+1}$. The following claim ensures such a choice of path $R_i$.
Let $S\defeq \{ v_j: j\in B^4( J_1\cup J_2\cup J'_2)\}$.

\begin{claim} \label{claim:finding_path}
For each $i\in I_3$, there exists a directed path $R_i = v^i_{j_1}v^i_{j_2}\ldots v^i_{j_{t}}$ in $T[W_i]- S$ such that $j_{k+1} - j_k\ge 100p$ for each $k\in [t-1]$, and $j_1 \le 100p$ and $j_t \geq n_i - 50p$. 
\end{claim}
\begin{claimproof} 
  Note that we have $|S|\le 20p$. Set $j_1$ to be the lowest integer $k$ such that $v^i_k \notin S$, then we have $j_1 \le 50p$. Suppose that $\ell \ge 1$ is the maximal integer such that there exist $j_1,j_2,\ldots,j_\ell$ satisfying that $v^i_{j_k} \notin S$ for all $k\in [\ell]$ and $j_{k+1} - j_k\ge 100p$ for all $k\in [\ell-1]$. 
    If $n_i - j_\ell > 100p$, then the property \ref{median_basic_2} of median order implies that there are at least $25p$ arcs $\arc{v^i_{j_\ell}v^i_k}$ with $j_\ell+ 50 p \le k \le n_i$ and at least one of them is not in $S$. This provides an arc $\arc{v^i_{j_{\ell}} v^{i}_{j_{\ell+1}}}$ contradicting the maximality of $\ell$, thus $j_\ell > n- 100p$. This proves the claim.
\end{claimproof}

Using this claim, we fix such a path $R_i$ for each $i\in I_3$. For each $i\in I_3$, let 
\(j(i) \defeq j_1\) where \(R_i = v_{j_1}\dots v_{j_{\ell}}\). Let
\[   J_3 \defeq \{ j(i) : i\in I_3\}
    \quad\text{and}\quad 
    A \defeq \{ v_{j}: j\in J_3\}.
\]
For each vertex $v_j \in \bigcup_{i\in I_3} R_i$, as $v_j\notin S$, we know that $v_{j-2}v_{j-1}v_jv_{j+1}v_{j+2}$ is a subpath of $P^2$. For each vertex $v_j\in (\bigcup_{i\in I_3} R_i) \setminus A$,
we apply \Cref{prop:removal_of_vertex} to replace the part $v_{j-2}v_{j-1}v_jv_{j+1}v_{j+2}$ with a shorter path $v_{j-2} y y' v_{j+2}$ with $\{y,y'\}=\{v_{j-1},v_{j+1}\}$. By repeating this, we obtain a path $P^*$ with $V(P^*)= V(P^2) \setminus \left(\bigcup_{i\in I_3}V(R_i) \setminus A\right)$.
Now, for each $i\in I_3$ and $j=j(i)\in J_3$, we replace the arc $\arc{v_{j}v_{j+1}}$ with the path $R_i R'_i v_{j+1}$. Note that $R_i$ starts from the vertex $v_i$.
Let $y_i$ be the last vertex of $R_i$ and $x'_i, y'_i$ be the first and last vertex of $R'_i$, respectively.
In fact, because the vertices in $R'_i$ belong to $Y_i$, $i\in I_3$ and $v_{j+1} \notin M$, \ref{U23} implies $v_{j+1} \in N^+_{T_{d_2(y'_i)}(y'_i)}$. By similar logic the vertex $y_i$ in $R_i$ lies in $N^-_{d_2(x'_i)}(x'_i)$. 
Using these facts, we color the arc $\arc{y_i x'_i}$ with the color $d_2(x'_i)$ and we color the arc $\arc{y'_i v_{j+1}}$ by the color $d_2(y'_i)$. 

 Repeating this for all $i\in I_3$ yield a desired partially colored directed path $P^3$ with $V(P^3) = W\cup U$ where all the uncolored arcs are in $T$. Moreover, all the colored arcs have distinct colors in \(\bigcup_{u\in U} D'_u\) as they are all colored using the colors  in \ref{U1} and \ref{U23}.
 
Let $I = I_1\cup I_2\cup I'_2\cup I_3$, then we know $|I|\leq 2p$.
Then, for all $i\notin I$, we know that $P^3[W_i]$ is $P_i$, which is the directed Hamilton path of $T[W_i]$ coming from the median order.

\paragraph{Step 3. Coloring the arcs in $W_i$ greedily and utilizing \Cref{lem:rainbowDHP_1st}.}

Since $T\subseteq T^{\alpha}_C[W]$, at least $\alpha n$ colors are available in each arc of $T$. 
Consider the set $E$ of arcs of the digraph $\bigl(\bigcup_{i\in I} P^3[\{w_{i-1},w_i\}\cup W_i]\bigr)$. 
As $|I|\leq 2p$, the set $E$ contains at most $2p\cdot (\mu n + 1) \leq \alpha n/4$ arcs. 
All the arcs in $E$ incident to a vertex in $U$ are already colored by using at most $3p$ colors.
The remaining arcs in $E$ belong to $T$, and thus, by greedily choosing a color for each arc, we can color the arcs of those parts using distinct colors. Let the resulting path be $\widetilde{P}$, where all parts except $\bigcup_{i\notin I} \widetilde{P}[\{w_{i-1},w_i\}\cup W_i]$ are colored. Let $C'\subseteq C$ be the set of the remaining colors not used in $\widetilde{P}$. Clearly, $|C'|\ge n+p-1 - |E| \ge (1-\alpha/4)n$ and $|C'|=\sum_{i\in [r]\setminus I} (|W_i|+1)$.
Recall that each arc $\arc{uv} \in T$ satisfies $\arc{uv}\in T_c$ for at least $\alpha n$ choices of $c\in C$, hence it satisfies $\arc{uv}\in T_c$ for at least $\alpha n - |E| -p \geq (\alpha/2) |C'|$ choices of $c\in C'$.

For each $i\notin I$, let $V_i= W_i\cup \{w_{i-1}, w_i\}$. 
We now take a random partition of $C'$ into $\{C'_i: i\in [r]\setminus I\}$ so that $|C'_i|=|W_i|+1$.
As $|C'_i| \geq \gamma \mu n$, a standard Chernoff bound on hypergeometric distributions together with a union bound over all arcs $\arc{uv}\in T$ yield that there exists a choice $\{C'_i : i\in [r]\setminus I\}$ satisfying the following.
\begin{equation}\label{eq: random color part}
    \text{each arc $\arc{uv}$ of $T$ belongs to $T_c$ for at least $(\alpha/4)|C'_i|$ choices of $c\in C'_i$.}
\end{equation}

Consider the collection of tournaments $\TT'_i:=\{T_c[V_i] : c\in C'_i\}$ on the vertex set $V_i$. 
Let $T^i = T[V_i]$. Note that \eqref{eq: random color part} implies that $T^i\subseteq \TT_i'^{\alpha/4}$.
Applying \Cref{lem:robust-H-partition} to $T^i-\{w_{i-1},w_i\}$, we obtain an $\Hpart(\mu|V_i|,\gamma)$-partition. As $w_{i-1}\outdir W_i$ and $W_i\outdir w_i$ in $T^i= T[V_i]$, we can apply \Cref{lem:rainbowDHP_1st} 
to this $\Hpart(\mu|V_i|,\gamma)$-partition with $w_{i-1}, w_i, \alpha/4$ playing the role of $w_0, w_{r}, \alpha$, respectively.
Then it yields a $C'_i$-rainbow path with the vertex set $\{w_{i-1},w_i\}\cup W_i$ from $w_{i-1}$ to $w_i$.
Repeating this for every $i\in [r]\setminus I$ and 
 concatenating these with the previously colored arcs yields a desired rainbow directed path. This finishes the proof of the lemma.
\end{proof}

\section{Rainbow brooms and DO-decompositions}
\label{sec:embedding_rainbow-brooms}

We now formally define the brooms mentioned in \Cref{sec:proofsketch}.

\begin{figure}[!htb]
    \centering
    \includegraphics{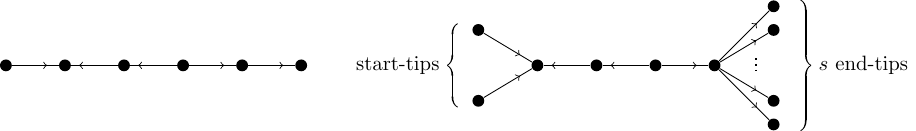}
    \caption{An example of a path \(P\) and a \((P,2,s)\)-broom.}
    \label{fig:broom}
\end{figure}

\begin{definition}[Broom]\label{def: broom}
For a path $P=x_1\dots x_{\ell}$ with $\ell \geq 2$ and integers $s_1,s_2$, a \emph{$(P,s_1,s_2)$-broom} is a digraph obtained from $P$ by blowing up $x_1$ into $s_1$ vertices and $x_\ell$ into $s_2$ vertices. When blowing up a vertex into $s$ vertices, the new arcs all keep the same orientation. See \Cref{fig:broom}.
If $\ell=1$, then we define $(P,s_1,s_2)$-broom be a set of $s_2$ vertices with no arcs.
We call the duplicates of $x_1$ \emph{start-tips} and the duplicates of $x_{\ell}$ \emph{end-tips} and a vertex is a \emph{tip} of the broom if it is either a start-tip or and end-tip.
If $s_1=1$ (or $s_2=1$) and there is only one start-tip (end-tip), then we call it the start-vertex (end-vertex).
If all start-tips belong to a vertex set $V_1$ and all end-tips belong to a vertex set $V_2$, then we say that the broom is \emph{from $V_1$ to $V_2$}.
The arcs adjacent to the duplicates of $x_1$ ($x_{\ell}$, respectively) are called \emph{start-tip-arcs} (\emph{end-tip-arcs}, respectively), and we call them together as \emph{tip-arcs}.
For a $(P,s_1,s_2)$-broom $F$, we let $\ell(F)\defeq \ell(P)=\ell-1$ be the length of the broom.
The path $x_2\dots x_{\ell-1}$ is called an \emph{internal path}. If the internal path of $F$ is disjoint from a vertex set $U$, we say that $F$ is internally disjoint from $U$.
\end{definition}

Note that if $P$ has length $1$, then the $(P,s_1,s_2)$-broom is a complete bipartite graph between $s_1$ vertices and $s_2$ vertices where all the arcs are consistently oriented as $P$. 

As we discussed in \Cref{sec:proofsketch}, we will find a broom and then choose one start/end-tip to obtain a path. In order to ensure that this path is rainbow, we may allow two incident tip-arcs in the broom to have the same color, but all other pairs of edges in the broom must have different colors. The following definition encapsulates this concept.

\begin{definition}[Near-rainbow coloring]
    A coloring $\varphi$ of a $(P,s_1,s_2)$-broom $F$ is \emph{near-rainbow} if all the start-tip-arcs have the same color and all the end-tip-arcs have the same color, and every path in $F$ from a start-tip to an end-tip is rainbow under $\varphi$. For two sets $C$ and $C'$ of colors, $\varphi$ is \emph{$(C,C')$-near-rainbow} if \(\varphi\) is near-rainbow and $C'\subseteq \varphi(F)\subseteq C$. If \(C'\) is empty, we simply say \(\varphi\) is \emph{\(C\)-near-rainbow}.
\end{definition}

Note that in a near-rainbow broom, all the paths from a start-tip to end-tip are rainbow paths using the same set of colors. The following proposition ensures the existence of brooms in a tournament.

   \begin{proposition}
    \label{clm:double-broom}
        For $\ell,s_1,s_2\in \mathbb{N}$ and \(n\ge \ell+2^{s_1+s_2}\), any $n$-vertex tournament $T$ contains a \(( \arc{P}_{\ell},s_1,s_2)\)-broom.
   \end{proposition}
    \begin{proof}
        If \(\ell=1\), consider a transitive subtournament $S$ in $T$ of size $s_1 + s_2$ with a median order \(\langle\omega_1,\dots,\omega_{s_1+s_2}\rangle\). Letting \(X\defeq\{\omega_1,\dots,\omega_{s_1}\}\) and \(Y\defeq\{\omega_{s_1+1},\dots,\omega_{s_1+s_2}\}\). As $X\outdir Y$ in $T$, we obtain a \((\arc{P}_{1},s_1,s_2)\)-broom. 
        
        Now suppose \(\ell\ge 2\). Consider a median order \(T=\langle\omega_1,\dots,\omega_n\rangle\). Let \(y\defeq\omega_{2s_1}\) and \(z\defeq\omega_{\ell+2s_1-2}\). 
As $2^{s_1+s_2}\geq 2s_1+2s_2$ for $s_1,s_2\geq 1$, we use \ref{median_basic_2} and \ref{median_basic_3} to
        choose \(Y\subseteq N^-(\omega_{2s_1},\{\omega_1,\dots,\omega_{2s_1-1}\})\) and \(Z\subseteq N^+(\omega_{\ell+2s_1-2},\{\omega_{\ell+2s_1-1},\dots,\omega_n\})\) of sizes \(s_1\) and \(s_2\), respectively. The path $\omega_{2s_1} \omega_{2s_1+1} \dots \omega_{\ell+2s_1-3}\omega_{\ell+2s_1-2}$ together with those sets $Y$ and $Z$
        gives an \((\arc{P}_{\ell},s_1,s_2)\)-broom with start-tips $Y$ and end-tips $Z$. 
    \end{proof}

Using the above proposition, we can obtain a large near-rainbow $(P,s_1,s_2)$-broom for each $P\in \{\arc{P}_{\ell}, \arcrev{P}_{\ell}\}$ in a given collection of tournaments if we have some additional vertices and colors. For our purpose, it is sufficient to only use it with $(s_1,s_2)=(2,50)$.
Throughout the rest of this paper, we fix a large constant 
\begin{equation}\label{eq: fix N0}
    N_0 \defeq 10 k_0 + 5000,
\end{equation}
where $k_0$ satisfies the following lemma.

\begin{lemma}
\label{lem:rainbow-dir-double-broom}
    There exists an integer \(k_0\) satisfying the following. Let \(\ell\ge 1\) and \(\TT\) be a collection of tournaments with
    \(|V(\TT)|,|\TT|\ge \ell+ k_0\). Then, \(\TT\) contains a near-rainbow $(P,2,50)$-broom for each $P\in \{\arc{P}_{\ell},\arcrev{P}_{\ell}\}$.
\end{lemma}
\begin{proof}
Choose $\mu, k_0>0$ so that $0< 1/k_0\ll \mu \ll 1$ holds.
We only prove the existence of a near-rainbow \((\arc{P}_{\ell},2,50)\)-broom as the other one follows by symmetry. Choose a tournament \(T\subseteq\TT^{1/2}\).
We start by finding a large broom in this uncolored tournament $T$.

If $\ell \le k_0 - 2^{120}$, then find a $(\arc{P}_{\ell},3,100)$-broom $F'$ in $T$ using \Cref{clm:double-broom}. 
Since at least half of the colors are available for each $100$ end-tip-arcs in $F'$, there exists a color $c_\fin$ available for at least fifty end-tip-arcs. We color fifty end-tip-arcs by the color $c_\fin$ and remove the other fifty end-tip-arcs.
Similarly, we color two of the start-tip-arcs in $F'$ with a color $c_\ini$ (different from $c_\fin$) and delete the other start-tip-arc.
For all remaining uncolored arcs, at least $\frac{1}{2}(\ell+k_0)-2 \geq \ell + 2^{119}-2> \ell$ colors are available, so we can greedily color them in a rainbow way to obtain a near-rainbow \((\arc{P}_{\ell},2,50)\)-broom.

Hence we may assume $\ell>  k_0- 2^{120}$, so we have $1/\ell \ll \mu$. 
By deleting vertices if necessary, we may assume that the number of vertices in $V(\TT)$ is exactly $\ell + 2^{102}+1$.

By \Cref{lem:robust-H-partition}, $T$ has a \(\Hpart(\mu \ell, 1/25)\)-partition \(\calH\defeq (W_0,\dots,W_{r+1},w_0,\dots,w_r)\). 
By \Cref{clm:double-broom}, there is a \((\arc{P}_{|W_0|-2^{101}},3,1)\)-broom \(L_0\) in \(T[W_0]\) and a \((\arc{P}_{|W_{r+1}|-2^{101}},1,100)\)-broom \(L_{r+1}\) in \(T[W_{r+1}]\). 
Let \(x_0\) be the end-tip of \(L_0\) and \(x_{r+1}\) be the start-tip of \(L_{r+1}\).

Similarly, as before, using the pigeonhole principle, we choose a color $c_\fin$ available for at least fifty end-tip-arcs of $L_{r+1}$. We color fifty end-tip-arcs by the color $c_\fin$ and remove the other fifty end-tip-arcs to obtain a \((\arc{P}_{|W_{r+1}|-2^{101}},1,50)\)-broom $L'_{r+1}$.
Again, we color two of the start-tip-arcs in $L_0$ with a color $c_\ini$ (different from~$c_\fin$) and delete the other start-tip-arc to obtain \((\arc{P}_{|W_{r+1}|-2^{101}},2,1)\)-broom $L'_0$.
Since at least $\frac{1}{2}(\ell+k_0)$ colors are available for the each uncolored arc in $L'_0\cup L'_{r+1} \cup \{\arc{x_0w_0}, \arc{w_{r}x_{r+1}}\} $ and $|L'_0|+|L'_{r+1}| \leq 3\mu \ell< \frac{1}{2}(\ell+k_0) -4$, we can find a rainbow coloring of the remaining arcs of $L'_0\cup L'_{r+1}\cup \{\arc{x_0w_0}, \arc{w_{r}x_{r+1}}\}$ avoiding the colors $c_\ini$ and~$c_\fin$. 

Let $C$ be the set of remaining unused colors.
The remaining \(\Hpart(\mu \ell, 1/25)\)-partition $(W_1,\dots,W_{r},w_1,\dots,w_{r-1})$ with two vertices $w_0, w_{r}$ satisfies $w_0\outdir W_1$ and $W_{r} \outdir w_{r}$ in $T$ and $T\subseteq \mathbb{T}_C[V_0]^{1/3}$ where $V_0=\bigcup_{i\in [r]}W_i \cup \{w_0,\dots, w_{r}\}$. 
As $1/\ell \ll \mu \ll 1$,  \Cref{lem:rainbowDHP_1st} yields a $C$-rainbow directed path from $w_0$ to $w_{r}$ covering all 
vertices in $V_0$. By concatenating this with the previously colored digraphs $L'_0$, $\arc{x_0w_0}$, $\arc{w_{r}x_{r+1}}$, and $L'_{r+1}$, we obtain a desired near-rainbow $(P,2,50)$-broom with the length exactly 
$\sum_{i\in [r]}(|W_i|+1) + \ell(L_0)+ \ell(L_1) = \ell+2^{102} - 2^{101}-2^{101} =\ell$. This proves the lemma.
\end{proof}

The following lemma obtains a large near-rainbow $(P,2,50)$-broom for each $P\in \{\arc{P}_{\ell}, \arcrev{P}_{\ell}\}$ with the additional property that the tips have large out-degree and in-degree outside of the broom.

\begin{lemma}\label{lem: rich directed path}
Let $0<1/\ell_0 \ll \varepsilon \ll 1/k \le 1$ and $\ell \in \mathbb{N}$ with $\ell \leq 10\ell_0$.
Suppose $\TT$ is a collection of tournaments with $|V|, |\Gamma|\geq \ell+ \ell_0$, where $V = V(\TT)$ and $\Gamma=\Gamma(\TT)$.
If a color set $D\subseteq \Gamma$ with $|D|\leq k$ is given, \(\TT\) contains a $(\Gamma\setminus D)$-near-rainbow $(P,2,50)$-broom $F$ for each $P\in \{\arc{P}_{\ell},\arcrev{P}_{\ell}\}$ such that each tip $v$ of $F$ satisfies $|N^{\sigma}_{T_d}(v, V\setminus V(F))| \geq \varepsilon \ell_0$ for all $d\in D$ and $\sigma\in \{+,-\}$.
\end{lemma}
\begin{proof}
As before, we assume $P=\arc{P}_{\ell}$ by symmetry.
Choose a tournament \(T\subseteq\TT^{1/2}_{\Gamma\setminus D}\).

We first consider the case $\ell\le \varepsilon \ell_0$. Let $V'\subseteq V$ be the set of vertices $v$ which has at least $3\varepsilon \ell_0$ out-neighbors and in-neighbors in each tournament $T_d$ with $d\in D$. 
\Cref{prop:num-of-small-indeg-vtxs_small} implies that 
\[
|V'| \geq |V| - |D| \cdot 5 \cdot 3\varepsilon\ell_0 \ge \ell_0 - 20k\varepsilon\ell_0> 4 \ell+1000,
\]
so we can use \Cref{prop:embed_tree} to find a copy of $(P,3,100)$-broom $F'$ in $\TT_{\Gamma\setminus D}^{1/2}[V']$. Similarly as before, using the pigeonhole principle, we choose a color $c_\fin$ available for at least fifty end-tip-arcs of $F'$ and a color $c_{\ini}$ (different from $c_\fin$) available for at least two start-tip-arcs of $F'$ and color them with those colors and delete the remaining tip-arcs. Using the definition of $T$ and the inequality that $\ell < k\varepsilon \ell_0 < \ell_0/100 < (\ell_0 - k)/ 2 \le |\Gamma\setminus D|/2$, it is easy to check that we can greedily give a $(\Gamma\setminus (D \cup \{c_\ini,c_\fin\}))$-rainbow coloring to the internal path of $F'$ to obtain a $(\Gamma\setminus D)$-near-rainbow $(P,2,50)$-broom $F$. 
As the tips of $F$ lies in $V'$, for every $d\in D$, we have
$$|N^{\sigma}_{T_d}(v, V\setminus V(F))| \geq 3\varepsilon \ell_0 - (\ell + 51) \geq \varepsilon \ell_0,$$
as desired.

Otherwise, we have $\ell> \varepsilon \ell_0$, thus $0<1/\ell \ll \varepsilon$ holds. Let $n\defeq |V|$.
We use \Cref{lem:robust-H-partition} to obtain a \(\Hpart( \varepsilon^{1/2} n, 1/25)\)-partition \(\calH\defeq (W_0,\dots,W_{r+1},w_0,\dots,w_r)\).
Thus, each $W_i$ has size between $\varepsilon^{1/2} n/25$ and $\varepsilon^{1/2} n$.

Using \Cref{prop:num-of-small-indeg-vtxs_small}, for each $i\in \{0,r+1\}$, we choose a set $W'_i\subseteq W_i$ of $(1- 500k\varepsilon^{1/2}) |W_i|\geq 100$ vertices $w$ that satisfy \begin{equation}\label{eq:minimum degree in W_i}
    |N^{\sigma}_{T_d}(w, W_i)| \geq 100\varepsilon^{1/2} |W_i|\geq 2 \varepsilon n
\end{equation} 
for all $\sigma\in \{+,-\}$ and $d\in D$. 

Note that $m\defeq |V\setminus (W_0\cup W_{r+1})|$ satisfies $m\geq n-2\varepsilon^{1/2} n \geq (1-2\varepsilon^{1/2})( \ell + \ell_0) >\ell$ as we have $\ell \leq 10\ell_0$.
For each $i\in [r]$, we take arbitrary subset $W'_i\subseteq W_i$ so that each set $W'_i$ has size 
either $\bigl\lfloor \frac{\ell-r-2}{r} \bigr\rfloor$ or $\bigl\lceil \frac{\ell-r-2}{r} \bigr\rceil$ so that
$\sum_{i\in [r]} (|W'_i|+1) = \ell-2$. 
This is possible as $m \geq \ell\geq r+2$. Moreover, since $\calH$ is  a \(\Hpart( \varepsilon^{1/2} n, 1/25)\)-partition, it is routine to check that  
\((W'_1,\dots, W'_{r}, w_1,\dots, w_{r-1})\) is a \(\Hpart( \varepsilon^{1/2} \ell, 1/100)\)-partition.

Choose three vertices $x$ in $W'_0$ and one hundred vertices $y$ in $W'_{r+1}$ and take the arcs of the form $x w_0$ and~$w_r y$. Similar to the previous case, using the definition of $T$ and the pigeonhole principle, we color two of the start-tip-arcs of the form $x w_0$ with a color $c_\ini \in \Gamma\setminus D$ and delete the other arcs from $W'_0$ to $w_0$, and then color fifty end-tip-arcs $w_r y$ by a color $c_\fin \in \Gamma\setminus (D\cup \{c_\ini\})$ and remove the other fifty arcs from $w_r$ to~$W'_{r+1}$.

Let $C$ be the set of remaining unused colors in $\Gamma \setminus D$, i.e., $C = \Gamma\setminus (D \cup \{c_\ini,c_\fin\})$.
As $(W'_1,\dots, W'_r,w_1,\dots, w_{r-1})$ is a  \(\Hpart(\varepsilon^{1/2} \ell, 1/100)\)-partition with two vertices $w_0, w_{r}$ satisfying $w_0\outdir W'_1$ and $W'_{r} \outdir w_{r}$ in $T$ and $T\subseteq \mathbb{T}_{C}[V'_0]^{1/3}$ where $V'_0=\bigcup_{i\in [r]}W'_i \cup \{w_0,\dots, w_{r}\}$. 
As $0<1/\ell \ll \varepsilon \ll  1$,  \Cref{lem:rainbowDHP_1st} yields a $C$-rainbow directed path from $w_0$ to $w_{r}$ covering all vertices in $V'_0$, yielding a $(\Gamma\setminus D)$-near-rainbow $(\arc{P}_{\ell},2,50)$-broom $F$.
Note that the tips of $F$ belong to $W'_{0}\cup W'_{r+1}$ and the vertices in $(W_{0}\cup W_{r+1}) \setminus (W'_{0}\cup W'_{r+1})$ are all outside $V(F)$. Hence, by \eqref{eq:minimum degree in W_i}, each tip $v$ of $F$ has at least $2\varepsilon n - 52 \geq \varepsilon n \ge \varepsilon \ell_0$ out-neighbors and in-neighbors in $V\setminus V(F)$ in each tournament $T_d$ for each $d\in D$. This proves the lemma.
\end{proof}

\subsection{DO-decomposition}

As we remarked in \Cref{sec:proofsketch}, we now define oscillating paths and the DO-decomposition.

\begin{definition}\label{def:osillating}
    Let \(P\) be a path. We say it is \emph{oscillating} if each of its blocks has length at most \(2\). Furthermore, an oscillating path is \emph{good} if \(\ell(P)\ge2\) and the last two arcs have opposite orientations.
\end{definition}

\begin{definition}\label{def:DO}
    For a path $P$, a decomposition \(P_1 P_2\dots  P_{2k} = P^\dirc_1 P^\oscl_1 \dots P^\dirc_k P^\oscl_k\) of $P$ into paths is a \emph{DO-decomposition (directed-oscillating decomposition)} if it satisfies the following. 
    \begin{proplist}
        \item \label{E1} Each path in the decomposition is nonempty possibly except the first one $P^\dirc_1$ and the last one $P^\oscl_k$.
        \item \label{E2} Each \(P^\dirc_i\) is a directed path and each $P^{\oscl}_i$ is an oscillating path for $i\in [k]$. Moreover, if $P^{\oscl}_i$ is nonempty, it is a non-directed path with at least two arcs.
        
        \item \label{cond:(lem:decomposing_path)_oscl-paths_shape} 
        For each $i\in [k-1]$, both \(P^\oscl_i\) and \(\rev(P^\oscl_i)\) are good oscillating paths. 
        If the last two arcs of $P$ have opposite orientations, then both \(P^\oscl_k\) and \(\rev(P^\oscl_k)\) are good.
        
        \item \label{cond:(lem:decomposing_path)_incident-paths} 
        For each $i\in [2,2k-2]$, the last arc of $P_i$ and the first arc of $P_{i+1}$ have the same orientation in $P$.
    \end{proplist}
\end{definition}

Indeed, the following lemma confirms that every path admits a DO-decomposition.

\begin{lemma}
\label{lem:decomposing_path}
    Every path \(P\) admits a DO-decomposition.
\end{lemma}

\begin{proof} 
If \(P\) has length at most two, then $P$ together with an empty path with an appropriate order is a DO-decomposition. Suppose that \(P\) has length at least three.

    Consider the blocks of $P$. Merging all consecutive blocks of length at most two, we obtain a decomposition $P= P_1Q_1P_2Q_2\dots P_{k} Q_k$ such that each $P_i$ is a block of length at least three (possibly except $P_1$ which may be an empty path) and each $Q_i$ is a (possibly empty) oscillating path. 
    For each $i\in [k]$, if $P_i$ is nonempty, let $P_{i}^{\dirc}$ be the path $P_i$ minus the first and last arc. 
    If $P_1$ is empty and the first two arcs of $Q_1$ have the same orientation, then let $P_1^{\dirc}$ be the first arc of $Q_1$. If $P_1$ is empty and the first two arcs of $Q_1$ have different orientations, then let $P_1^{\dirc}$ be the empty path.
    
    Let $P_{i}^{\oscl}$ be
    the path between the last vertex of $P_i^{\dirc}$ and the first vertex of $P_{i+1}^{\dirc}$. It is clear that each $P_{i}^{\oscl}$ is an oscillating path, as it is obtained by attaching one arc of opposite orientation at one or both ends of an oscillating path $Q_i$. 
    Moreover, \ref{E1} and \ref{E2} are clear, and \ref{cond:(lem:decomposing_path)_incident-paths} is straightforward from our choice.
    By our construction, the first two arcs in $P_1^{\oscl}$ have different orientations, and for each $i\in [2,k]$, the first two arcs in $P_{i}^{\oscl}$ belong to the consecutive blocks $P_i$ and $Q_i$, thus have different orientations.
    Similarly, for each $i\in [k-1]$, the last two arcs in $P_{i}^{\oscl}$ have different orientations. Hence $P_{i}^{\oscl}$ and $\rev({P_{i}^{\oscl}})$ are both good oscillating paths.
    If, in addition, the last two arcs of $P$ have different orientations, then they both belong to $P_{k}^{\oscl}$, hence the last two arcs in $P_{k}^{\oscl}$ also have different orientations, hence \ref{cond:(lem:decomposing_path)_oscl-paths_shape} follows. This proves the lemma.
\end{proof}

\subsection{Near-rainbow brooms}

In this subsection, we show that, for an oscillating path $P$, finding a rainbow copy of $P$ is possible even with some requirements on the choices of colors. In fact, as we discussed in \Cref{sec:proofsketch}, we want to find a near-rainbow $(P,2,50)$-broom instead of the path $P$.
The following two propositions find near-rainbow brooms for very short oscillating paths $P$.

\begin{proposition}
\label{lem:rainbow-short-free-oscl-double-broom}
Let \(\TT\) be a collection of tournaments and \(V_1,V_2\subseteq V(\TT)\) be disjoint sets with \(|V_1|\ge 50\) and and \(|V_2|\ge 300\). 
For any oscillating path $P$ of length three and a  set \(B\subseteq\Gamma(\TT)\) of three colors, \(\TT[V_1\cup V_2]\) contains a $B$-near-rainbow $(P,2,50)$-broom from $V_1$ to $V_2$.
\end{proposition}
\begin{proof}
Let $(\sigma_1,\sigma_2, \sigma_3)$ be the orientations of three arcs of $P$ and let $B=\{1,2,3\}.$ By \Cref{prop:num-of-small-indeg-vtxs_small}, find two vertices $u_1, u_2\in V_1$ such that \(d_{T_1}^{-\sigma_1}(u_i,V_1)\ge 5\) for \(i\in [2]\) and find a set $W$ of $50$ vertices $w$ in $V_2$ such that \(d_{T_3}^{\sigma_3}(w,V_2)\ge 51\).
Without loss of generality, assume $(u_1u_2)^{\sigma_2}\in T_3$.

If there are $i\in [2]$ and $w\in W$ such that $(u_iw)^{\sigma_2}\in T_2$, then taking two vertices in $N^{-\sigma_1}_{T_1}(u_i)$ and $50$ vertices in $N^{\sigma_3}_{T_3}(w,V_2)$, we obtain a desired $B$-near-rainbow $(P,2,50)$-broom. Thus we assume that $(u_iw)^{-\sigma_2}\in T_2$ for all $i\in[2]$ and $w\in W$.

As $P$ is an oscillating path, we have $\sigma_1\neq \sigma_2$ or $\sigma_2\neq \sigma_3$.
If $\sigma_1\neq \sigma_2$, then choose an arc $(ww')^{\sigma_2}\in T_1[W]$, then taking arcs $(u_1w)^{\sigma_1}=(u_1w)^{-\sigma_2}\in T_2$ and $(u_2w)^{\sigma_1}\in T_2$ and $50$ vertices in $N_{T_3}^{\sigma_3}(w')\setminus \{w\}$ yields a desired broom.
If $\sigma_2\neq \sigma_3$, then the arc $(u_1u_2)^{\sigma_2}\in T_3$ with two vertices in $N_{T_1}^{-\sigma_1}(u_1,V_1)\setminus\{u_2\}$ and $50$ vertices in $W\subseteq N_{T_2}^{-\sigma_2}(u_2) =N_{T_2}^{\sigma_3}(u_2)$, yields a desired broom.
\end{proof}

\begin{proposition}
\label{lem:rainbow-short-lim-oscl-double-broom}
Let \(\TT\) be a collection of tournaments and  \(V_1,\{v\}\subseteq V(\TT)\) be disjoint subsets with \(|V_1|\ge 50\).
    Let \(P= x_1\dots x_{\ell+1}\) be a good oscillating path of length $2\leq \ell \leq 4$, let $P' \defeq P-x_{\ell+1}$, and let \(B\subseteq\Gamma(\TT)\) be a color set of size \(\ell\) and $b\in B$.
    Then \(\TT[V_1\cup\{v\}]\) contains one of the following.
    \begin{enumerate}
        \item A \(B\)-near-rainbow \((P,2,1)\)-broom from $V_1$ to $\{v\}$, or
        \item a \((B\setminus\{b\})\)-near-rainbow \((P',2,1)\)-broom from $V_1$ to $\{v\}$.
    \end{enumerate}
\end{proposition}
\begin{proof}
As $P$ is a good oscillating path, without loss of generality, we may assume $\arc{x_{\ell-1}x_{\ell}}, \arcrev{x_{\ell}x_{\ell+1}} \in P$.
Let $\sigma\in \{+,-\}$ with $(x_1x_2)^{\sigma}\in P$ and let $P'' \defeq P'-x_1 = x_2\dots x_\ell$ and $P''' \defeq P''-x_\ell = x_2\dots x_{\ell -1}$.
We further assume that $B=[\ell]$ and $b=\ell-1$.
Take $U\defeq N_{T_\ell}^{+}(v,V_1)$ and $U'\defeq N_{T_\ell}^{-}(v,V_1)$.

If $|U|\geq 25$, then \Cref{prop:num-of-small-indeg-vtxs_small} yields a set $W\subseteq U$ of five vertices $w$ satisfying $d_{T_1}^{-\sigma}(w,U) \geq 5$. Apply \Cref{prop:simple-embed-rainbow-path} to obtain a $[2,\ell-1]$-rainbow path $P''=w_2\dots w_{\ell}$ in $W$. Adding two vertices in $N_{T_1}^{-\sigma}(w_2)\setminus V(P'')$ and the vertex $v \in N_{T_\ell}^{-}(w_{\ell})$, we obtain a desired $[\ell]$-near-rainbow $(P,2,1)$-broom from $V_1$ to $v$.

Suppose $|U'|\geq 25$.
If $\ell=2$, then \(v\) together with two vertices in \(U'\) yields a desired \(\{2\}\)-near-rainbow \((P',2,1)\)-broom to \(v\). Assuming $\ell\geq 3$, apply \Cref{prop:num-of-small-indeg-vtxs_small} to obtain a set $W\subseteq  U'$ of five vertices $w$ satisfying $d_{T_1}^{-\sigma}(w,U) \geq 5$.
Again, apply \Cref{prop:simple-embed-rainbow-path} to obtain a $[2,\ell-2]$-rainbow path $P'''=w_2\dots w_{\ell -1}$ in $W$. Adding two vertices in $N_{T_1}^{-\sigma}(w_2)\setminus V(P''')$ and the vertex $v \in N_{T_\ell}^{+}(w_{\ell -1})$, we obtain a desired ($[\ell]\setminus\{\ell-1\}$)-rainbow $(P',2,1)$-broom to $v$.
\end{proof}

Repeatedly using these two propositions, 
we can prove the following lemma. It states that, for an oscillating path $P$, we can find a $P$-broom even with some restrictions of colors and with specified start vertex set and end vertex set. The proof is straightforward, so we only provide a sketch.

\begin{lemma}
\label{lem:rainbow-oscl-path/broom}
    Let \(\TT\) be a collection of tournaments with \(|V(\TT)|\ge \ell(P)+1000\) and $V_1,V_2 \subseteq V(\TT)$ be disjoint vertex sets with $|V_1|\geq 50$. 
    Let \(P\) be a non-empty oscillating path and $B\subseteq \Gamma(\TT)$ be a set of $\ell(P)$ colors. 
    \begin{proplist}
        \item \label{case:(lem:rainbow-oscl-path/broom)_loscl-end}
        If $P$ is a good oscillating path and $|V_2|\geq 2$, then \(\TT\) has a \(B\)-near-rainbow \((P,2,1)\)-broom  from \(V_1\) to~\(V_2\). 

        \item \label{case:(lem:rainbow-oscl-path/broom)_foscl-end}
        If \(|V_2|\ge 300\) and $\ell(P)$ is divisible by three, then \(\TT\) has a \(B\)-near-rainbow \((P,2,50)\)-broom from $V_1$ to~$V_2$.
    \end{proplist}
\end{lemma}
\begin{proof}[Proof sketch]
Assume $|V_1| = 50$ and $|V_2|\leq 300$ by deleting vertices from $V_1$ and $V_2$ if necessary.
We decompose $P=P_1\dots P_{k}$ into subpaths $P_1,\dots, P_{k-1}$ of length three and $P_k$ of length between $2$ and $4$. 
Also, partition the color set $B$ into $B_1,\dots, B_k$ with $|B_i|=\ell(P_i)$ for each \(i\).
Let $F_0$ be the vertex set $V_1$, which is a $(P_0,2,50)$-broom where $P_0$ is the empty path. 
We repeatedly apply \Cref{lem:rainbow-short-free-oscl-double-broom} to obtain, for each $i\leq k-1$, a $B_i$-near-rainbow $(P_i,2,50)$-broom $F_i$ in $\TT_{B_i}\setminus V_2$ whose start-tips lie in the set of end-tips of the previous broom $F_{i-1}$ and it does not intersect with the previous brooms except that its start-tips belong to the set of end-tips of $F_{i-1}$. For the last step, we apply \Cref{lem:rainbow-short-lim-oscl-double-broom} for \ref{case:(lem:rainbow-oscl-path/broom)_loscl-end} to find a $(P_k,2,1)$-broom $F_k$ or apply \Cref{lem:rainbow-short-free-oscl-double-broom} for \ref{case:(lem:rainbow-oscl-path/broom)_foscl-end} to find a $(P_k,2,50)$-broom $F_k$, where $F_k$ is $B_k$-near-rainbow broom from the end tips of $F_{k-1}$ to $V_2$. Concatenating those brooms $F_i$ and removing redundant vertices yields a desired $B$-near-rainbow broom from $V_1$ to $V_2$.
\end{proof}

Using this lemma and the DO-decomposition, we can find a near-rainbow $P$-broom for a general path $P$ as shown in the following lemma. Some of its assumptions may seem technical, as this lemma will be applied frequently in various settings. Recall that $N_0$ is the constant we fixed in~\eqref{eq: fix N0}.

\begin{lemma}
\label{lem:rainbow_broom_Deb_version}
    Let \(\TT\) be a collection of tournaments and $V_2 \subseteq V(\TT)$. 
    Let $n\defeq |V(\TT)|$ and $\ell$ be such that $n\geq \ell+ \frac{1}{4}N_0$. 
    Let \(P\) be a nonempty path of length \(\ell\) admitting a DO-decomposition \(P^\dirc_1 P^\oscl_1 \dots P^\dirc_k P^\oscl_k\). 
    Let $V_1\subseteq V(\TT)$ be a vertex set disjoint with $V_2$ such that if $\ell(P^\dirc_1) = 0$, then $|V_1|\geq 50$ and otherwise, $|V_1|\ge \ell(P^\dirc_1) + \frac{1}{10}N_0$. 
    Let \(B,C\subseteq \Gamma(\TT)\) be color sets satisfying \(B\subseteq C\), \(|C|\ge \ell+ \frac{1}{8}N_0\), and \(|B|\le \sum_{i\in [k]} \ell(P^\oscl_i)\). 
    \begin{proplist}
        \item \label{case:(lem:rainbow_broom_Deb_version)_loscl-end}
        If \(P^\oscl_k\) is good oscillating and \(|V_2|\ge 2\), then \(\TT\) has a \((C,B)\)-near-rainbow \((P,2,1)\)-broom \(F\) from $V_1$ to~$V_2$.
        
        \item \label{case:(lem:rainbow_broom_Deb_version)_dir-end}
        If \(P^\oscl_k\) is an oscillating path with $3\mid \ell(P^{\oscl}_k)$, then \(\TT\) has a \((C,B)\)-near-rainbow \((P,2,50)\)-broom $F$ from $V_1$. 
        Moreover, if \(P^\oscl_k\) is nonempty and $|V_2|\ge 300$, then we can further ensure that $F$ is from $V_1$ to $V_2$.

        \item \label{new G2} 
        Let $\varepsilon, \ell_0$ be such that $0<1/\ell_0 \ll \varepsilon \ll 1/N_0$ and let $D\subseteq \Gamma(\TT)\setminus C$ be a set of at most $N_0$ colors.
        If $|C|\geq \ell+\ell_0$ and $n\geq \ell+ \ell_0$ holds in addition, and $P$ satisfies \(\ell(P^\oscl_k)= 0\) and \(\ell(P^{\dirc}_k)\leq 10\ell_0\),
        then \(\TT\) has a \((C,B)\)-near-rainbow \((P,2,50)\)-broom $F$ from $V_1$ where each end-tip $v$ of $F$ satisfies
        $|N^{\sigma}_{T_d}(v, V(\TT)\setminus V(F))| \geq \varepsilon \ell_0$ for each $d\in D$ and $\sigma\in \{+,-\}$.
    \end{proplist}
\end{lemma}

\begin{proof}[Proof of \Cref{lem:rainbow_broom_Deb_version}]
As before, by deleting some vertices if necessary, we may assume that $|V_2|\leq 300$.
Also, by adding some colors in $C$ to $B$ if necessary, we may assume that 
$|B|=\sum_{i\in [k]} \ell(P^{\oscl}_i)$.
We arbitrarily partition $B$ into $B_1,\dots, B_{k}$ so that $|B_i|=\ell(P^{\oscl}_i)$ for each $i\in [k]$. 

For $i\in [k]$, let  $P^i=P_1^{\dirc} P_1^{\oscl} P_2^{\dirc} \dots P_{i}^{(d)}$ be the subpath of $P$ obtained by concatenating the first $2i-1$ paths of the given DO-decomposition of $P$. In what follows, we will prove \ref{case:(lem:rainbow_broom_Deb_version)_loscl-end}--\ref{new G2} simultaneously including the moreover part, and we will sometimes have the same notations denoting objects with different definitions based on which one among \ref{case:(lem:rainbow_broom_Deb_version)_loscl-end}--\ref{new G2} is being considered.
We will inductively show the following. 
\begin{itemize}
    \item For every $i \in [k]$, there is a $(C,\bigcup_{j\in [i-1]}B_{j})$-near-rainbow $(P^i,2,50)$-broom $F_i$ from $V_1$ in $\TT\setminus V_2$ without using any colors in $\bigcup_{j\in [i,k]}B_{j}$. 
\end{itemize}
We start by proving it for $i=1$. 
If $P^1=P_1^{\dirc}$ is an empty path, then any arbitrary set of $50$ vertices in $V_1$ is a $(P^1,2,50)$-broom without any edges, so a desired $F_1$ exists.
If $P^1=P_1^{\dirc}$ is non-empty, then using the facts 
\(|V_1|\ge \ell(P^\dirc_1)+N_0/10\) and 
\(|C\setminus B| \ge \ell(P^\dirc_1)+N_0/10\), 
apply \Cref{lem:rainbow-dir-double-broom} to find a $C\setminus B$-near-rainbow
\((P^\dirc_1,2,50)\)-broom \(F_1\) in \(V_1\).

For the inductive step, suppose we have $i_*\in [k-1]$ and a $(C,\bigcup_{j\in [i_*-1]}B_{j})$-near-rainbow $(P^{i_*},2,50)$-broom $F_{i_*}$ from $V_1$ in $\TT\setminus V_2$ without using any colors in $\bigcup_{j\in [i_*,k]}B_{j}$. 
Let $Z$ be the set of end-tips of $F_{i_*}$ and let $\Phi_{i_*}$ be the set of colors used in $F_{i_*}$.
As $i_*<k$, using the facts 
\(|V(\TT)\setminus (V(F_{i_*})\cup V_2)|\ge \ell(P^\dirc_{i_*+1})+N_0/5\) and 
\(|C\setminus (\Phi_{i_*}\cup B)| \ge \ell(P^\dirc_{i_*+1})+N_0/10\), 
apply \Cref{lem:rainbow-dir-double-broom} to find a $C\setminus (\Phi_{i_*}\cup B)$-near-rainbow
\((P^\dirc_{i_*+1},2,50)\)-broom \(F^*\) in \(\TT \setminus (V(F_{i_*})\cup V_2)\). Let $Y^*$ be the set of start-tips of $F^*$. 

If the conditions in \ref{new G2} are met and $i_*=k-1$, then we have \(|V(\TT)\setminus (V(F_{i_*})\cup V_2)|\ge \ell(P^\dirc_{i_*+1})+\ell_0/2\) and 
\(|C\setminus (\Phi_{i_*}\cup B)| \ge \ell(P^\dirc_{i_*+1})+\ell_0/2\), 
and thus apply \Cref{lem: rich directed path} with $2\varepsilon$ in place of $\varepsilon$ to find a 
$C\setminus (\Phi_{i_*}\cup B)$-near-rainbow
\((P^\dirc_{i_*+1},2,50)\)-broom \(F'\) in \(\TT \setminus (V(F_{i_*})\cup V_2)\), then every tip of $F'$ has at least $2\varepsilon (n-|V(F_{i_*}) \cup V_2|) \geq \varepsilon \ell_0$ out-neighbors and in-neighbors in $V'_k= V(\TT)\setminus (V(F_{i_*})\cup V_2)$.
Let $Y'$ be the set of start-tips of $F'$.

As $i_*<k$, by the definition of the DO-decomposition, \(P^\oscl_{i_*}\) is a nonempty good oscillating path with $|B_{i_*}|=\ell(P^{\oscl}_{i_*})$. 
If the conditions in \ref{new G2} are not met or $i_*<k-1$, let $W=\emptyset$. 
If the conditions in \ref{new G2} are met and $i_*=k-1$, then for each of fifty two tips of $F'$, we choose (not necessarily disjoint) $\varepsilon \ell_0$ out-neighbors and $\varepsilon \ell_0$ in-neighbors in each of $T_d$ for $d\in D$ to obtain a vertex set $W$ of size at most $52|D|\varepsilon \ell_0 \leq 52 N_0 \varepsilon \ell_0 < \ell_0 - N_0/2$.

In either case, as there are at least $|P_i^{\oscl}|+ N_0/5 \geq |P_i^{\oscl}| + 1000 $ vertices outside $V(F_{i_*})\cup V(F^*)\cup V_2\cup W$ (or $V(F_{i_*})\cup V(F')\cup V_2\cup W$), we can apply \Cref{lem:rainbow-oscl-path/broom}-\ref{case:(lem:rainbow-oscl-path/broom)_loscl-end} to the path $P_i^{\oscl}$ to obtain a $B_{i_*}$-near-rainbow $(P_i^{\oscl},2,1)$-broom $F''$ in $\TT$ from $Z$ to $Y^*$ (or $Y'$) which is internally disjoint from $V(F_{i_*})\cup V(F^*)\cup V_2\cup W$ (or $V(F_{i_*})\cup V(F')\cup V_2\cup W$).
Concatenating $F, F'', F^*$ in order (or $F, F'', F'$ in order) and deleting redundant vertices yields a $(C,\bigcup_{j\in [i_*]}B_{j})$-near-rainbow $(P^{i_*+1},2,50)$-broom $F_{i^* +1}$ from $V_1$ in $\TT\setminus V_2$ without using any colors in $\bigcup_{j\in [i_*+1,k]}B_{j}$.  
As this procedure can be done as long as $i_*<k$, this yields a $(C,\bigcup_{i\in [k-1]}B_{i})$-near-rainbow $(P^k,2,50)$-broom $F_k$ from $V_1$ in $\TT\setminus V_2$ without using any colors in $B_k$, as desired. 
Let $Z^*$ be the set of end-tips of $F_k$.
If $P^{\oscl}_k$ is empty, then the first sentence of \ref{case:(lem:rainbow_broom_Deb_version)_dir-end} holds. If the conditions in \ref{new G2} are met, then the definition of $W$ naturally yields the conclusion in \ref{new G2}.

To prove \ref{case:(lem:rainbow_broom_Deb_version)_loscl-end} and \ref{case:(lem:rainbow_broom_Deb_version)_dir-end} for the case when $P^{\oscl}_k$ is not empty,
we apply \Cref{lem:rainbow-oscl-path/broom}-\ref{case:(lem:rainbow-oscl-path/broom)_loscl-end} and \ref{case:(lem:rainbow-oscl-path/broom)_foscl-end} to $\TT\setminus (V(F_k)\setminus Z^*)$ to find a $B_{k}$-near-rainbow $(P^{\oscl}_k,2,1)$ broom and a $B_{k}$-near-rainbow $(P^{\oscl}_k,2,50)$-broom from $Z^*$ to $V_2$ in each case, respectively. This is possible because $\TT\setminus (V(F_k)\setminus Z^*)$ has at least $\ell(P_k^{\oscl})+ N_0/5$ vertices.
Moreover, this ensures that they are not intersecting $F_{k}$.
Again concatenating this new broom with $F_k$ and deleting redundant arcs incident to $Z^*$, we obtain desired near-rainbow brooms, confirming \ref{case:(lem:rainbow_broom_Deb_version)_loscl-end} and \ref{case:(lem:rainbow_broom_Deb_version)_dir-end}. 
\end{proof}

\section{Rainbow paths with only short blocks}
\label{sec:shortblocks}

We now prove \Cref{thm:rainbow-rosenfeld_path_short-blocks}. We recall the statement for convenience.

\ShortBlocks*

\begin{proof}[Proof of \Cref{thm:rainbow-rosenfeld_path_short-blocks}]

\begin{figure}[htb]
    \centering
    \includegraphics{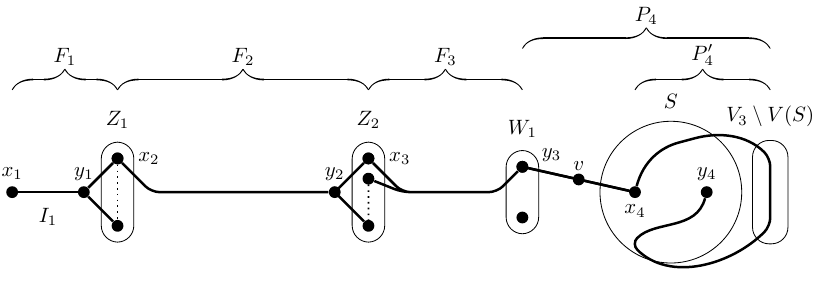}
    \caption{A description for the proof of \Cref{thm:rainbow-rosenfeld_path_short-blocks}.}
    \label{fig:thm@rainbow-rosenfeld_path_short-blocks)-pf_description}
\end{figure}

We choose $\beta,\gamma>0$ satisfying the hierarchy
\begin{equation*}\label{eq:(thm:rainbow-rosenfeld_path_short-blocks)-pf param hierarchy}
    0<1/n \ll \gamma \ll \beta \ll 1.
\end{equation*}
We start by finding a decomposition \(P_1P_2P_3P_4\) of \(P\) into nonempty paths satisfying the following and some length requirements that will be mentioned later. For each \(i\in \{2,3\}\), there is a DO-decomposition \(P_i=P_{i,1}^\dirc P_{i,2}^\oscl \cdots P_{i,k_i}^\dirc P_{i,k_i}^\oscl\) which satisfies the following.
\begin{proplist}
    \item \label{P2P_J} $P_{2,1}^{\dirc}$ is empty and $3\divides \ell(P_{2,k_2}^{\oscl})$;
    \item \label{P3P_J} $P_{3,1}^{\dirc}$ is empty and $P_{3,k_3}^{\oscl}$ is a good oscillating path.
\end{proplist}
To find such a decomposition \(P_1P_2P_3P_4\), 
let $P= v_1\dots v_\ell$ = 
$P_1^\dirc P_1^\oscl \dots P_s^\dirc P_s^\oscl$ be a DO-decomposition of $P$.
For each $i\in [s-1]$, \ref{E1} and \ref{E2} implies that
$P_{i}^{\oscl}= u^i_1\dots u^i_{f_i}$ satisfies $f_i\geq 3$. 
Let 
\[
    A_i \defeq \{ u^i_1, u^i_4, \dots, u^i_{3\lceil (f_i -2)/3 \rceil -2}\}
    \quad\text{and}\quad
    A \defeq \bigcup_{i\in [s-1]} A_i.
\]

Note that for every $i\in [s-1]$, we have 
\begin{equation}\label{eq:number of arcs left at the end}
2\le f_i - (3\lceil (f_i -2)/3 \rceil -2) \le 4.
\end{equation}
Enumerate $A$ as $v_{i_1},\dots, v_{i_a}$ with $i_1<\dots <i_a$. As each block has length at most $(\log n)^{1/2}$, we have $i_{j+1}-i_j \leq (\log n)^{1/2}+4$ for each $j\in [a-1]$.
Choose the maximum $a_1$ and the maximum $a_2>a_1$ satisfying the following
\[
    i_{a_1} \le \beta n+1
    \quad\text{and}\quad i_{a_2}-i_{a_1} \le (1-11\beta+\gamma)n.
\]
Let $i'$ be the maximum number less than $n - \frac{1}{3}\log n$ such that the two arcs $v_{i'-2}v_{i'-1}$ and $v_{i'-1} v_{i'}$ have the opposite orientations. Then we have $i' \geq n-\frac{1}{3}\log n - (\log n)^{1/2} -4 \geq n-\frac{1}{2}\log n$.
Let 
\[
    P_1 = v_1\dots v_{i_{a_1}},
    \quad P_2 = v_{i_{a_1}}\dots v_{i_{a_2}},
    \quad P_3 = v_{i_{a_2}}\dots v_{i'},
    \quad\text{and}\quad
    P_4 = v_{i'}\dots v_{\ell}.
\]
Note that the DO-decomposition $P_1^\dirc P_1^\oscl \dots P_s^\dirc P_s^\oscl$ naturally yields a DO-decomposition for each of $P_2$  and $P_3$ by using the subpaths $P_j^\dirc \cap P_i$ and $P_j^{\oscl}\cap P_i$ whenever they contain at least one vertex. By our choices of $a_1,a_2,a_3$, 
it is straightforward that these DO-decomposition of $P_2$ and $P_3$ satisfy \ref{P2P_J} and \ref{P3P_J}. 
In particular, the lower bound in \eqref{eq:number of arcs left at the end} ensures that $P_{2,1}^{\dirc}$ and $P_{3,1}^{\dirc}$ are empty. 
The choice of the vertices in $A$ ensures that $3\divides \ell(P_{2,k_2}^{\oscl})$ and the choice of $a_3$ ensures that the last two arcs of $P_3$ are opposite and thus \ref{cond:(lem:decomposing_path)_oscl-paths_shape} and \ref{cond:(lem:decomposing_path)_incident-paths} imply that $P_{3,k_3}^{\oscl}$ is a good oscillating path.
Finally, we slightly adjust the values of \(\beta\) and \(\gamma\) if necessary to assume the following while maintaining the hierarchy $1/n\ll\gamma \ll \beta \ll 1$.
\begin{gather*}
    \ell(P_1)=\beta n +1,
    \quad \ell(P_2) = (1-11\beta +\gamma) n - N_0 -4,
    \\ \quad 9\beta n\leq \ell(P_3)\leq  (10\beta - \gamma/2 ) n,
    \quad \text{and} \quad \frac{1}{3}\log n \leq \ell(P_4) \leq \frac{1}{2}\log n.
\end{gather*}

\paragraph{Preparation of a color absorber and a transitive tournament.} 
We fix three colors $1,2,n$ and a tournament 
\(T\subseteq\TT_{[3,n-1]}^{1/2}\). We plan to use the colors $1$ and $2$ for the first two arcs of $P_4$.
Using \Cref{prop:num-of-small-indeg-vtxs_small}, take a vertex set $V_1$ consisting of $5\beta n$ vertices $v$ with $d^+_{T_n}(v) \geq n/4$.

By applying \Cref{prop:embed_tree} on $T[V_1]$, we get a \((P_1,1,100)\)-broom  in \(T[V_1]\).
Using the pigeonhole principle, choose a color $d \in [3,n-1]$ such that at least $50$ end-tip-arcs are in $T_d$. Let $Z_1$ be a set of such $50$ end-tips and let $F_1$ be the \((P_1,1,50)\)-broom obtained by deleting the end-tips not in $Z_1$. Let $x_1$ be the start vertex of $F_1$ and let $y_1\in V(F_1)$ be the vertex that corresponds to the second-to-last vertex of \(P_1\). 

Let $I_1$ be the subpath of $F_1$ from $x_1$ to $y_1$.
We construct a color absorber using this path $I_1$. By \Cref{lem:absorber}, there are disjoint sets \(A,C\subseteq[3,n-1]\setminus\{d\}\) with \(|A|=(\beta-\gamma)n\) and \(|C|=10\beta n\) such that
\begin{equation}\label{eq: color absorb short block}
    \text{for any \(C'\subseteq C\) of size \(\gamma n\), there is an (\(A\cup C'\))-rainbow coloring of \(I_1\).}
\end{equation}
Let 
$B\defeq [3,n-1] \setminus (A\cup C\cup \{d\}),$
then we have 
\begin{equation}\label{eq: Bsize}
    |B| = (1-11\beta +\gamma)n -4 = \ell(P_2)+ N_0.
\end{equation}

As we have prepared a color absorber, we now prepare a subtournament for the transitive absorption. 
Let $\sigma_1, \sigma_2$ be the orientation of the first two arcs of $P_4$.

Let \(V_2\defeq N_{T_n}^+(x_{1})\setminus V(F_1)\). Since $x_1 \in V_1$, we have \(|V_2|\ge n/4 - (\beta n + 51)\ge n/5\). Thus, \Cref{prop:num-of-small-indeg-vtxs_small} yields a vertex \(v\in V_2\) such that it has at least $n/20$ in-neighbors and at least $n/20$ out-neighbors in $V_2$.

Choose a set \(W_1\subseteq N^{-\sigma_1}_{T_1}(v,V_2)\) of size \(2\).
Consider the vertex set $U= N^{\sigma_2}_{T_2}(v,V_2)\setminus W_1$ and the tournament $T'=\TT[U]^{1/2}_{C}$. As $T'$ has at least $n/20-2$ vertices, it contains a transitive subtournament of size at least $\log(n/20-2)\geq \frac{1}{2}\log n \geq \ell(P_4)-N_0$.
Let $S$ be a transitive subtournament of $T'$ on exactly $\ell(P_4)-N_0$ vertices. 
As each arc in $S$ belongs to $T'$, there are at least $\frac{1}{2}|C|\geq 5\beta n$ colors available for each arc. Hence, we can greedily color all arcs in $S$ using at most $(\log n)^2$
colors to obtain a $C$-rainbow coloring $\varphi_1$ of $S$. Let $W\defeq W_1\cup \{v\} \cup V(S)$. We have $|W|= \ell(P_4)-N_0+3$.

\paragraph{Finding a near-rainbow \((P_2,1,50)\)-broom $F_2$ using most of the colors in $B$.} 
By \ref{P2P_J}, we can apply \Cref{lem:rainbow_broom_Deb_version}-\ref{case:(lem:rainbow_broom_Deb_version)_dir-end} to obtain a $B$-near-rainbow $(P_2,1,50)$-broom $F_2$ from $Z_1$ in the vertex set $V(\TT)\setminus(W\cup V(F_1))$. This is possible by \eqref{eq: Bsize} as $V(\TT)\setminus(W\cup V(F_1))$ has more than $n-\beta n-\log n > \ell(P_2)+ N_0$ vertices.
Let $x_2\in Z_1$ be the start-tip and $Z_2$ be the set of end-tips of $F_2$ and let $\varphi_2$ be the $B$-near-rainbow coloring.
Let 
$$B_1 = B\setminus \varphi_2(F_2) \enspace \text{and} \enspace 
D= (C \setminus \varphi_1(S))\cup B_1.$$ 
We have $|B_1| \leq N_0$ and $|D|\ge |C \setminus \varphi_1(S)| > 10\beta n -(\log n)^2 > \ell(P_3)+ N_0$.

\paragraph{Finding a near-rainbow \((P_3,2,1)\)-broom $F_3$ exhausting the colors in $B$.} 
Let $U= V(\TT)\setminus (V(F_1)\cup V(F_2)\cup W)$.
We apply \Cref{lem:rainbow_broom_Deb_version}-\ref{case:(lem:rainbow_broom_Deb_version)_loscl-end} to $\TT[U]$ to obtain a $(D,B_1)$-near-rainbow $(P_3,2,1)$-broom $F_3$ from $Z_2$ to a vertex $y_3\in W_1$.
Indeed, we can check that such an application is possible.
First, $P_{3,1}^{\dirc}$ is empty and $P_{3,k_3}^{\oscl}$ is a good oscillating path by \ref{P3P_J}.
Second, we have enough vertices and colors as $V(\TT)\setminus (V(F_1)\cup V(F_2)\cup W)$ has $n-\ell(P_1)-\ell(P_2) -|W|- 99 = n- \ell(P_1)-\ell(P_2)-\ell(P_4)+N_0-102 = \ell(P_3)+N_0 - 101$ vertices
and $|D|> \ell(P_3)+ N_0$. 
Lastly, as every block in $P$ has length at most $(\log n)^{1/2}$, we have $k_3\geq \frac{\ell(P_3)}{(\log n)^{1/2}} \geq \frac{9\beta n}{(\log n)^{1/2}} > N_0  \geq  |B_1|$. 
Since \(P^\oscl_{3,j}\) is nonempty for every \(j\in [k_3]\), we have \(\sum_{j=1}^{k_3} \ell(P^\oscl_{3,j})\geq k_3\ge |B_1|\). Thus, the application of \ref{case:(lem:rainbow_broom_Deb_version)_loscl-end} is possible, and we obtain $F_3$.

Now, we concatenate $F_1,F_2,F_3$ and discard the redundant vertices in $Z_1\cup Z_2$ to obtain a path $Q$. Recall that all the end-tip-arcs of $F_1$ are colored by the same color $d$, and $F_2, F_3$ are near-rainbow, so this removal of the redundant vertices does not change the set of colors used to partially color $F_1\cup F_2\cup F_3$. Thus the path $Q$ has the first uncolored segment $I_1$, and the rest is colored using all colors in $B \cup\{d\}$ and some additional colors in $C$.

\paragraph{Attaching a rainbow $P_4$ at the end of $Q$ using the vertex-absorbing transitive tournament $S$.}
Let $V_3\defeq V\setminus (V(Q)\cup \{v\})$, then this set has size $n-\ell(P_1)-\ell(P_2)-\ell(P_3)-2 = \ell(P'_4)+1 = |W|+ N_0-4 = |V(S)|+N_0-1$, where $P'_4$ is the path obtained from $P_4$ by deleting first two arcs. 

The transitive tournament $S$ lies in $V_3$ and has size $\ell(P_4)-N_0$, so \(|V_3\setminus V(S)| =  N_0-1\) holds. 
Let $C^*$ denote the set of all colors in $C$ that are not yet used in coloring $Q$ or $S$. Then \(|C^*|\ge 10\beta n-((\log n)^2+(10\beta-\gamma/2)n) > (\gamma/3)n\). 
Let \(T'\) be a tournament on \(V_3\) whose arc set is a subset of 
\[E(S)\cup 
\{\arc{uv}\in \TT[V_3]_{C^*}^{1/2} : u\notin V(S)\text{ or }v\notin V(S)\}.\] 
Since \(|E(T')| <  (\log n)^2 < |C^*|/2\), by a greedy choice we can give \(E(T')\setminus E(S)\) a \(C^*\)-rainbow coloring \(\varphi'\).

Since $\ell(P_4)\ge \frac{1}{3}\log n > 2^{N_0+3}$ and  
\(b(P_4)\ge \frac{\frac{1}{3}\log n}{(\log n)^{1/2}}\ge 5(N_0-1)\) and \(|V(T')\setminus V(S)|= N_0-1 \), we can apply \Cref{lem:Ham-path_p-excep_start-end-in-S} to find a copy of $P'_4$ in $T'$ which starts and ends at $S$ while containing all vertices in $V_3\setminus V(S)$.
As the arcs in $T'$ are already colored, this naturally gives rise to the coloring of $P'_4$.
Let \(x_{4}\) and \(y_{4}\) be the start- and end-vertex of this copy of $P'_4$, respectively. 
Then using our construction of $W_1$ and $S$,
we add two arcs $(y_3 v)^{\sigma_1}$ from $T_1$ and $(vx_4)^{\sigma_2}$ from $T_2$. As these two arcs together with $P'_4$ yield a copy of $P_4$, this together with $Q$ yields a copy of $P$ where every arc except the ones in $I_1$ are colored using colors in $[n-1]\setminus A$ such that the set $C'$ of all unused colors in $[n-1]$ is a subset of $C$.
Moreover, as the last vertex of $P$ is in $S \subseteq N^{+}_{T_n}(x_1)$, the last vertex of $P$ is an out-neighbor of the first vertex of $P$ in $T_n$.
Finally, we can apply \eqref{eq: color absorb short block} to color arcs in $I_1$ using colors in $A\cup C'$ to finish the coloring of $P$. This finishes the proof of \Cref{thm:rainbow-rosenfeld_path_short-blocks}.
\end{proof}

\section{Rainbow paths with long blocks}
\label{sec:rainbow-paths_long-blocks}
We now prove \Cref{thm:rainbow-rosenfeld-path_long-block-end}. We recall its statement below.
We remind that for a given path $P$, we defined $\shift(P)$ as the path obtained from $P$ by removing the last arc and appending it at the start of the path in the reverse direction.

\LongBlock*

\begin{proof}[Proof of \Cref{thm:rainbow-rosenfeld-path_long-block-end}]

\begin{figure}[htb]
    \centering
    \includegraphics{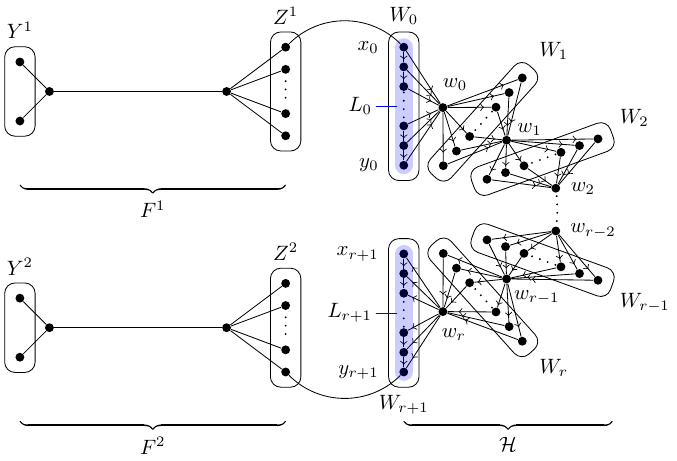}
    \caption{A description for the proof of \Cref{thm:rainbow-rosenfeld-path_long-block-end}.}
\end{figure}

Let $V\defeq V(\TT)$ and \(\Gamma\defeq \Gamma(\TT)\). By considering $\rev(P)$ instead of $P$ if necessary, we may assume that there is a longest block which is not the first block in $P$.
Let $Q^*$ be the last longest block. By flipping the orientations of all arcs in all tournaments in $\TT$ and all arcs in $P$ if necessary, we may assume that this block $Q^*$ is forwardly oriented.
Consider the decomposition $P=P_{1} Q P'_{2}$ where $Q$ is $Q^*$ minus the first arc and the last arc of $Q^*$. If $Q^*$ is the last block, then we do not delete the last arc so that $P'_{2}$ is the empty path in this case.
Let $P_2 \defeq \rev(P'_{2})$. 
Since $P$ is not a directed path, $P_1$ is not an empty path, and letting $\ell\defeq \ell(Q)$, we have $(\log n)^{1/2}-2\leq \ell \leq n-3$.
We plan to find disjoint near-rainbow brooms of subpaths of \(P_1\) and \(P_2\), and connect their end-tips using a directed path obtained from $\Hpart$-absorption. 
For $i\in [2]$, if $P_{i}$ is nonempty, let 
\(P^\dirc_{i,1} P^\oscl_{i,1} \dots P^\dirc_{i,k_i} P^\oscl_{i,k_i}\) be a DO-decomposition of \(P_{i}\). Since $\ell(P^\dirc_{1,1}) + \ell \le n-1$ and $\ell(P^\dirc_{1,1})\le \ell(Q^*)\le \ell -1$, we have $\ell(P^\dirc_{1,1})\le n/2$.

Choose $\varepsilon, \mu, \gamma$ so that $0<1/n \ll \mu \ll \varepsilon, \gamma \ll 1/N_0$. Choose a set $C_0\subseteq \Gamma$ of $N_0$ colors. By \Cref{prop:num-of-small-indeg-vtxs_small}, there are at least $n-4\varepsilon^{1/2} N_0 n > 5n/6\ge \ell(P^\dirc_{1,1}) + \ell/3$ vertices having both out-degree and in-degree at least $\varepsilon^{1/2} n$ in $T_{c}$ for all $c\in C_0$. Fix such is a set $U'$ consisting of exactly $\ell(P^\dirc_{1,1}) + \ell/3$ such vertices. 
We choose a vertex set $W'\subseteq V$ of size exactly $\ell/5$ satisfying the following. 
\newlength{\richseteqlen}
\settowidth{\richseteqlen}{Every vertex $u\in U'$ has at least $\varepsilon \ell$ out-neighbors and in-neighbors in $W'$ in every $T_c$ with $c\in C_0$.}
\begin{equation}\label{eq: rich set}
    \parbox{\richseteqlen}{Every vertex $u\in U'$ has at least $\varepsilon \ell$ out-neighbors and in-neighbors in $W'$ in every $T_c$ with $c\in C_0$.}
\end{equation}
Indeed, if we include every vertex of $V$ in $W'$ independently with probability $\ell/(10n)$, then by the standard use of Chernoff bound and union bound, we have that with positive probability, the number of vertices in $W'$ is at most $\ell/5$ and $W'$ satisfies \eqref{eq: rich set}. Choose one such $W'\subseteq V$. Now, by adding arbitrary vertices to $W'$ (note that this does not affect the property in \eqref{eq: rich set}), we can ensure that $W'$ has size exactly $\ell/5$.

Choose arbitrary three disjoint sets \(U^1_\ini, U^1_\fin, U^2_\fin \in U'\setminus W'\) such that $|U^1_\ini| = \ell(P^\dirc_{1,1}) + N_0/10$ and $|U^i_\fin| = 300$ for $i\in [2]$. Indeed, this is possible since $|U'\setminus W'|\ge \ell(P^\dirc_{1,1}) + \ell/3 - \ell/5 \ge \ell(P^\dirc_{1,1}) + N_0/10 + 600$. Let $U\defeq U^1_\ini \cup U^1_\fin \cup U^2_\fin$.

\paragraph{Finding near-rainbow brooms for the subpaths of \(P_1\) and \(P_2\).}
Choose two disjoint sets $V_1, V_2$ in $V\setminus W'$ $V_1, V_2$ and choose disjoint color sets \(\Psi_1,\Psi_2\subseteq \Gamma\setminus C_0\) satisfying the following.
\[U^1_{\ini}\cup U^1_{\fin}\subseteq V_1,
\quad U^2_{\fin}\subseteq V_2,
\quad |V_i| = \ell(P_i)+ \ell/5,
\quad \text{and}
\quad |\Psi_i|=\ell(P_i)+ \ell/5 \quad \text{for each } i\in [2].\]
Define $U^2_{\ini} \defeq V_2\setminus U^2_{\fin}$. Note that $|U^2_{\ini}|\ge \ell(P^\dirc_{2,1}) + N_0/10$. 
We do the following process for each \(i\in[2]\).
If $P_i$ is nonempty, as the last two arcs of $P_{i}$ have opposite orientations, $P^\oscl_{i,k_i}$ is a nonempty good oscillating path by \ref{cond:(lem:decomposing_path)_oscl-paths_shape}. We further decompose $P^{\oscl}_{i,k_i}$ into $R_{i,1}R_{i,2}$ such that $3\mid \ell(R_{i,1})$ and $\ell(R_{i,2})\in \{2,3,4\}$.
Let $R_{i}^* \defeq P^\dirc_{i,1} P^\oscl_{i,1} \dots P^\dirc_{i,k_i} R_{i,1}$, then \(P_i = R_i^* R_{i,2}\). 
Note that $P^\dirc_{i,1} P^\oscl_{i,1} \dots P^\dirc_{i,k_i} R_{i,1}$ is a DO-decomposition of $R_{i}^*$ with $R_{i,1}$ being the last oscilating path.
If $R_{i,1}$ is nonempty, apply \Cref{lem:rainbow_broom_Deb_version}-\ref{case:(lem:rainbow_broom_Deb_version)_dir-end} to obtain a \(\Psi_i\)-near-rainbow \((R_i^*,2,50)\)-broom \(F_i\) in \(\TT[V_i]\) from \(U^i_\ini\) to \(U^i_\fin\). 

If $P_{i}$ is nonempty but $R_{i,1}$ is an empty path, then we apply \Cref{lem:rainbow_broom_Deb_version}-\ref{new G2} with $\TT_{\Psi_i\cup C_0}[V_i],\ell(R_i^*), \ell/5, C_0$ playing the roles of $\TT,\ell, \ell_0, D$, respectively. 
Then, we obtain a \(\Psi_i\)-near-rainbow \((R_i^*,2,50)\)-broom \(F_i\) in \(\TT[V_i]\) from \(U^i_\ini\) where each end-tip of $F_i$ has at least $\varepsilon \ell/5$ out-neighbors and in-neighbors in $V_i\setminus V(F_i)$ for each tournament $T_c$ with $c\in C$. Indeed, such an application is possible as the last directed path $P^{\dirc}_k$ has length at most $\ell \leq 10 \cdot (\ell/5)$ from the definition of $\ell$. 

When $i=2$, in both of the above two cases based on whether $R_{2,1}$ is empty or not, we extract an \(\Psi_2\)-near-rainbow \((R_2^*,1,50)\)-broom from $F_2$ and redefine $F_2$ to be this new broom.

At this point, for each $i\in [2]$, if $P_i$ is nonempty, let \(Y_i\subseteq U^i_\ini\) and \(Z_i\) be the set of start-tips and end-tips of \(F_i\), respectively, and let \(\Phi_i\subseteq \Psi_i\) be the set of colors used for \(F_i\). Note that $|Y_1| = 2$, $|Y_2| = 1$, and $|Z_1| = |Z_2| = 50$.

If $P_{i}$ is an empty path, which is possible only when \(i=2\), we define the empty objects for convenience: let \(Y_2=Z_2\defeq \emptyset\), let \(R_2^*\) and \(R_{2,2}\) be empty paths, and \(F_2\) be an empty broom with no vertex or edge. Also, let \(\Phi_2\) be an empty color set.

From the above process, we now have the near-rainbow brooms $F_1$ and $F_2$ which are vertex-disjoint with \(V(F_1), V(F_2)\subseteq V\setminus W'\), and using disjoint sets of colors \(\Phi_1, \Phi_2\subseteq \Gamma(\TT)\setminus C_0\). 
Let $V'\defeq V\setminus V(F_1\cup F_2)$ and $\Phi'\defeq \Gamma\setminus (C_0\cup \Phi_1\cup \Phi_2)$. Then we have $|\Phi'|= n-1-\sum_{i\in [2]}(\ell(P_i)-\ell(R_{i,2})) - N_0 = \ell + \ell(R_{1,2})+\ell(R_{2,2}) - N_0$.
Whether $R_{i,1}$ is empty or not, the following holds.
\begin{align}\label{eq: rich}
\begin{split}
    &\text{Every vertex in $Y_1\cup Z_1\cup Z_2$ (i.e., start- and end-tips of $F_1$ and end-tips of $F_2$)}\\
&\text{have at least $\varepsilon \ell/5$ out-neighbors and in-neighbors in $V'$ in every tournament $T_c$ with $c\in C_0$.}
\end{split}
\end{align}
Indeed, if $R_{i,1}$ is nonempty, then the application of \Cref{lem:rainbow_broom_Deb_version}-\ref{case:(lem:rainbow_broom_Deb_version)_dir-end} ensures \(Z_i\subseteq U^i_\fin\) and thus \eqref{eq: rich set} implies \eqref{eq: rich}. On the other hand, if $R_{i,1}$ is nonempty, \eqref{eq: rich} is simply a property obtained from the application of \Cref{lem:rainbow_broom_Deb_version}-\ref{new G2}.

\paragraph{Construction of a good $\Hpart$-partition.}
We now construct an $\Hpart$-absorber. 
Fix a tournament \(T\subseteq \TT[V']_{\Phi'}^{1/2}\) and let \(\calH = (W_0,\dots,W_{r+1},w_0,\dots,w_r)\) be a good \(\Hpart(\mu\ell, \gamma)\)-partition of \(T\) obtained by \Cref{cor:robust-H-partition-avoid-first-and-last-parts}, which exists since \(|V'|=n-|V(F_1)|-|V(F_2)| \ge n - (\ell(P_1)+51)-(\ell(P_2)+50)\ge \ell-200\).
By taking a median order, we choose spanning directed paths \(L_0\) and \(L_{r+1}\) in \(T[W_0]\) and \(T[W_{r+1}]\), respectively. 
For $i\in \{0,r+1\}$, let \(x_{i}\) and \(y_{i}\) be the start-vertex and the end-vertex of \(L_i\), respectively. 
Arbitrarily color the arcs in $L\defeq L_0\cup L_{r+1}\cup \{ \arc{y_0w_0},\arc{w_r x_{r+1}}\}$ using colors from $\Phi'$ so that $L$ is rainbow. This is possible as $L$ contains less than $2\mu \ell +2 \leq \frac{1}{2}|\Phi'|$ arcs.
Let $\Phi''$ be the set of all colors not used in $F_1\cup F_2\cup L$.

\paragraph{Connecting $F_1\cup F_2$ with $L$ using rainbow $R_{1,2}$ and $R_{2,2}$.}
Choose arbitrary disjoint sets $D_1, D_2\subseteq \Phi''$ of size $\ell(R_{1,2})$ and $\ell(R_{2,2})$, respectively.
Let $u_1 = x_0$ and $u_{2} = y_{r+1}$.

For each $i\in [2]$, we do the following.
If \(R_{i,2}\) is empty (which only happens when $i=2$ and $P_2$ is empty), then let \(F'_i\) be an empty broom with no vertex or edge.
Otherwise, if $R_{i,2}$ is nonempty, we apply \Cref{lem:rainbow-short-lim-oscl-double-broom} to the path $R_{i,2}$ to obtain a $D_i$-near-rainbow $(R'_{i},2,1)$-broom $F'_{i}$ in $\TT[ Z_{i}\cup \{ u_i \}]$ from $Z_i$ to the vertex $u_i$, where $R'_{i}$ is either $R_{i,2}$ or the path obtained from $R_{i,2}$ by deleting the last arc. Note that the construction of our decomposition $P_1QP'_2$ yields that the deleted arc belongs to the long block $Q^*$.
Hence, it has the forward orientation same as $Q$ if $i=1$ and the backward orientation same as $\rev(Q)$ if $i=2$.

The above process yields near-rainbow brooms \(F'_1\) and \(F'_2\). We concatenate \(F_1\), \(F'_{1}\), \(L_0\), and \(\arc{y_{0}w_{0}}\), and remove redundant vertices in $Z_1$ to obtain a near-rainbow $(R_1^* R'_{1} \arc{P}_{\ell_1}, 2,1)$-broom $F_{1}^*$ from $Y_1$ to $w_0$, where \(\ell_1\defeq \ell(L_0)+1\).
Similarly, we concatenate \(F_{2}\), \(F'_{2}\), \(\rev(L_{r+1})\), and \(\arcrev{x_{r+1} w_r}\), and remove redundant vertices in $Z_2$ to obtain a near-rainbow $(R_2^* R'_{2} \arcrev{P}_{\ell_2},1,1)$-broom (or a \((\arcrev{P}_{\ell_2},1,1)\)-broom if \(P_2\) is empty) $I_{2}$ from \(Y_2\) to \(w_r\), where \(\ell_2\defeq \ell(L_{r+1})+1\). Note that $I_{2}$ is a rainbow path, and $\rev(I_2)$ is a concatenation of some forwardly directed path with $P'_2$. 

Let $y_1,y_2$ be the start-tips of $F_{1}^{\ast}$ and choose a color $c\in \Phi''$ not used in $F_{1}^{\ast}\cup I_2$.
Without loss of generality, assume $\arcrev{y_1y_2}\in T_c$.
Let $I^1_1$ be the rainbow path obtained from $F_{1}^{\ast}$ by deleting $y_2$, and let $I^2_1$ be the rainbow path obtained by adding the arc $\arcrev{y_1y_2}$ colored with $c$ at the beginning of the path obtained from $F_{1}^{\ast}$ by deleting $y_1$.
Both paths $I^1_1$ and $I^2_1$ start at the same vertex $y_1$ and use the same colors except that $I^2_1$ additionally uses the color $c$. The path $I^1_1$ is isomorphic to $P_1$ concatenated with some directed path at the ends, whereas $I^2_1$ is isomorphic to a concatenation of a backward arc and $I^1_1$. 
Note that the last vertex of $I_1^i$ is $w_0$, and the last vertex of $I_2$ is $w_r$.

Note that for each $i\in [2]$, the vertices outside $V'\cup V(I_1^{i})\cup V(I_2)$ are the vertices in \(Y_1\cup Z_1\cup Z_2\) that are removed from the start- and end-tips of $F_1$ and end-tips of $F_2$. Now, we will use $\Hpart$-absorption for these vertices using the property in \eqref{eq: rich}.

\paragraph{$\Hpart$-absorption.}
Let $\Phi^*$ be the set of all colors not used in $I^1_1\cup I_2$. Then, $\Phi^*$ contains all the colors in $C_0$ and the previously defined tournament $T$ satisfies $T\subseteq \TT[V']^{1/4}_{\Phi^*}$ because $|\Phi^*|\geq \ell(Q) - |W_0|-|W_{r+1}|-2\geq \ell -3\mu \ell$ and $\frac{1}{2}|\Phi'| - (|\Phi'|-|\Phi^*|) \geq \ell - 3\mu \ell - \frac{1}{2}(\ell + \ell(R_{1,2})+ \ell(R_{2,2}) - N_0)
\geq (\frac{1}{2}-4\mu)\ell \geq \frac{1}{4}|\Phi'| \geq \frac{1}{4}|\Phi^*|$.
Let 
\[
    U^* \defeq (Y_1\cup Z_1\cup Z_2)\setminus V(I^1_1 \cup I_2)
    \quad\text{and}\quad
    W \defeq V'\setminus (W_0\cup W_{r+1}).
\]
Since \(W'\subseteq V'=W\cup W_0\cup W_{r+1}\), at least $|W'|-2\mu \ell$ vertices of $W'$ lies in $W$.
Also, \(|U^*|\le |Y_1\cup Z_1\cup Z_2|\le 102\), and \eqref{eq: rich} implies that each vertex in $U^*$ has at least $\frac{1}{5} \varepsilon \ell - 2\mu \ell > \frac{1}{6}\varepsilon \ell$ in-neighbors and out-neighbors in $W$. 
In addition, $(W_1,\dots, W_r, w_1,\dots, w_{r-1})$ is a robust $\Hpart$-partition and $w_0\outdir W_1$ and $W_r\outdir w_{r+1}$. Moreover,
\begin{align*}
    |\Phi^*| &= |\Gamma|-\ell(I^1_1)-\ell(I_2) = (n+1)-|V(I^1_1\cup I_2)| \quad \text{and} \\
    |U^*| &= |V(F_1\cup F_2)|+|W_0\cup W_{r+1}|+|\{w_0,w_r\}|-|V(I^1_1\cup I_2)| \\
        &= n-|V'|-|V(I^1_1\cup I_2)|+|W_0\cup W_{r+1}|+2 \\
        &=|\Phi^*| +1 - (|V'|-|W_0\cup W_{r+1}|),
\end{align*}
so \(|\Phi^*|=|U^*|+|W|-1\).
Hence, we can apply \Cref{cor:rainbowDHP_digphs-attached} with the following parameters
\statementinput
    {U^*, W, \Phi^*, C_0, T, \calH, w_{0}, w_r, \varepsilon/6, 1/4}
    {U, W, C, D_u, T, \calH, w_0, w_r, \varepsilon, \alpha}
to get a rainbow directed Hamilton path \(J\) of \(\TT[V']\) from $w_0$ to $w_{r}$ using exactly the colors in $\Phi^*$.

Concatenating $I^1_1$, $J$, and $\rev(I_2)$, 
we obtain a rainbow spanning path whose left and right ends are isomorphic to $P_1$ and $P'_2$, and the middle part is a forwardly oriented direct path. Since the whole path has length exactly $n-1$, the final path is isomorphic to $P$.

If $P_2$ is an empty path, then we do the final $\Hpart$-absorption again with $I_1^2$ instead of $I^1_1$, then we obtain a rainbow copy of $\arcrev{e} P_1 Q'$ where $Q'$ is $Q$ minus the last arc and $\arcrev{e}$ is a backwardly oriented edge. Moreover, it has the same end vertices as the rainbow copy of $P$ obtained above. Thus, it yields the desired rainbow copy of $\shift(P)$ sharing the same end vertices with $P$. This finishes the proof of \Cref{thm:rainbow-rosenfeld-path_long-block-end}.
\end{proof}

\section{Rainbow Hamilton cycles}
\label{sec:cycle}

Here we prove \Cref{thm:rainbow-rosenfeld-cycle}. As mentioned in \Cref{sec:proofsketch}, this follows from \Cref{thm:rainbow-rosenfeld_path_short-blocks,thm:rainbow-rosenfeld-path_long-block-end}, but we need to separately treat the case where the cycle is a directed path plus one arc with the opposite direction.

\begin{proof}[Proof of \Cref{thm:rainbow-rosenfeld-cycle}]
    Let \(V\defeq V(\TT)\) and \(\Gamma(\TT)=[n]\). Let \(\calC\) be an arbitrary oriented Hamilton cycle that is not consistently oriented.  We divide into three cases depending on the length \(\ell\) of a longest block of \(\calC\). 

    \begin{Cases}
        \item \(\ell< (\log n)^{1/2}\).
        Decompose \(\calC\) into a path \(P\) of length \(n-1\) and a forward arc from the start-vertex of $P$ to the end-vertex of $P$. Apply \Cref{thm:rainbow-rosenfeld_path_short-blocks} to get an \([n-1]\)-rainbow path \(Q\) isomorphic to \(P\) which starts from \(u\) and ends at \(v\) such that \(\arc{uv}\in T_n\). This gives a transversal of \(\calC\) in $\TT$.

        \item \((\log n)^{1/2}\le \ell < n-1\).
        Decompose \(\calC\) into a path \(P\) of length \(n-1\) whose last block of \(P\) has length \(\ell\), and an arc from the start-vertex of $P$ to the end-vertex of $P$ with the same orientation as the last arc of $P$. Note that \(P\) is not a directed path since \(\ell<n-1\). By \Cref{thm:rainbow-rosenfeld-path_long-block-end}, we can find \(u,v\in V\) and \([n-1]\)-rainbow paths \(Q_1\) isomorphic to \(P\) and \(Q_2\) isomorphic to \(\shift(P)\) such that both \(Q_1\) and \(Q_2\) start from \(u\) and end at \(v\). If \(\arc{uv}\in T_n\), then this arc colored with $n$ together with \(Q_1\) is a transversal of \(\calC\); otherwise, the arc $\arcrev{uv}$ colored with $n$ together with \(Q_2\) is a transversal of \(\calC\), as desired.

        \item \(\ell=n-1\).
        Choose $\mu,\gamma>0$ such that \(0<\frac{1}{n}\ll \mu\ll \gamma\ll 1\) holds. Fix a tournament \(T\subseteq \TT_{[n-1]}^{1/2}\), and let \(\calH=(W_0,\dots,W_{r+1},w_0,\dots,w_r)\) be an \(\Hpart(\mu n, \gamma)\)-partition. 
        Apply \Cref{prop:modify_hamilton_path} on $T[W_0]$ to obtain two Hamilton paths $P_0$ and $P'_0$ starting with the same vertex $x_{0}$ such that $P_0$ is a directed path and every arc in $P'_0$ except exactly one has the same orientation as $P_0$. Let \(L_0\) be the path which is \(P_0\) if \(\arc{x_{0} y_{r+1}}\in T_n\) and $P'_0$ otherwise. Let $y_0$ be the last vertex of $L_0$.
        Moreover, let \(L_{r+1}\) be a directed Hamilton path in $T[W_{r+1}]$ from a vertex $x_{r+1}$ to a vertex $y_{r+1}$.
        
        We greedily assign pairwise distinct colors from $[n-1]$ to the arcs in $L_0\cup L_{r+1}\cup \{\arc{y_0w_0}, \arc{w_{r}x_{r+1}} \}$. Let $C$ be the set of used colors and let $W\defeq V\setminus (W_0\cup W_{r+1})$.
        We apply \Cref{lem:rainbowDHP_1st} to find a rainbow path from $w_0$ to $w_{r}$ on vertex set $W$ using exactly the colors in $[n-1]\setminus C$.
        Concatenating this with $L_0$, $L_{r+1}$, and the arc between $x_{0}$ and $y_{r+1}$ colored with $n$ gives us rainbow copy of \(\calC\). This proves \Cref{thm:rainbow-rosenfeld-cycle}.
        \qedhere
    \end{Cases}
\end{proof}

\section{Concluding remarks}\label{sec:concluding remarks}
Finally, we mention a few open questions and future directions for research. 

\paragraph{Exceptional tournaments} As mentioned in the introduction, Havet and Thomass\'{e}~\cite{havet2000oriented} and Zein~\cite{zein2022oriented} determined the list of all exceptional tournaments not containing every orientation of Hamilton path and not containing every orientation of Hamilton cycle except possibly the directed one, respectively. 
In a similar vein, it will be interesting to find all exceptional collections of tournaments in \Cref{thm:rainbow-rosenfeld-path,thm:rainbow-rosenfeld-cycle}.

\paragraph{Complexity questions} Various algorithmic aspects have been considered to find oriented Hamilton paths in tournaments. For example, a polynomial algorithm is found in \cite{bang-jensen1992polynomial} to decide if a directed Hamilton path exists between two specified vertices in a tournament. A similar result is obtained in \cite{hell2002antidirected} for the Hamilton paths with an orientation where every two adjacent edges have opposite directions. 
When the end vertices are not specified, it is shown in 
\cite{hell1983complexity} that the complexity of finding a Hamilton path in tournaments behaves differently depending on the orientation of the path. It will be interesting to extend such results to the setting of rainbow~paths.

\printbibliography

@preamble{ " \newcommand{\noop}[1]{} " }

@article {grunbaum1971antidirected,
    AUTHOR = {Gr\"{u}nbaum, Branko},
     TITLE = {Antidirected {H}amiltonian paths in tournaments},
   JOURNAL = {J. Combinatorial Theory Ser. B},
  FJOURNAL = {Journal of Combinatorial Theory. Series B},
    VOLUME = {11},
      YEAR = {1971},
     PAGES = {249--257},
      ISSN = {0095-8956},
   MRCLASS = {05C20},
  MRNUMBER = {291022},
MRREVIEWER = {R.\ Cori},
       DOI = {10.1016/0095-8956(71)90035-9},
       URL = {https://doi-org-ssl.access.yonsei.ac.kr/10.1016/0095-8956(71)90035-9},
}

@article {rosenfeld1972antidirected,
    AUTHOR = {Rosenfeld, Moshe},
     TITLE = {Antidirected {H}amiltonian paths in tournaments},
   JOURNAL = {J. Combinatorial Theory Ser. B},
  FJOURNAL = {Journal of Combinatorial Theory. Series B},
    VOLUME = {12},
      YEAR = {1972},
     PAGES = {93--99},
      ISSN = {0095-8956},
   MRCLASS = {05.60},
  MRNUMBER = {285452},
MRREVIEWER = {A. C. Shershin},
       DOI = {10.1016/0095-8956(72)90035-4},
       URL = {https://www.sciencedirect.com/science/article/pii/0095895672900354},
}

@article {thomassen1973antidirected,
    AUTHOR = {Thomassen, Carsten},
     TITLE = {Antidirected {H}amilton circuits and paths in tournaments},
   JOURNAL = {Math. Ann.},
  FJOURNAL = {Mathematische Annalen},
    VOLUME = {201},
      YEAR = {1973},
     PAGES = {231--238},
      ISSN = {0025-5831},
   MRCLASS = {05C20},
  MRNUMBER = {349467},
MRREVIEWER = {L. S. Mel'nikov},
       DOI = {10.1007/BF01427945},
       URL = {http://link.springer.com/10.1007/BF01427945},
}

@article {rosenfeld1974antidirected,
    AUTHOR = {Rosenfeld, Moshe},
     TITLE = {Antidirected {H}amiltonian circuits in tournaments},
   JOURNAL = {J. Combinatorial Theory Ser. B},
  FJOURNAL = {Journal of Combinatorial Theory. Series B},
    VOLUME = {16},
      YEAR = {1974},
     PAGES = {234--242},
      ISSN = {0095-8956},
   MRCLASS = {05C35},
  MRNUMBER = {340101},
MRREVIEWER = {R. L. Hemminger},
       DOI = {10.1016/0095-8956(74)90069-0},
       URL = {https://www.sciencedirect.com/science/article/pii/0095895674900690},
}

@incollection {petrovic1984antidirected,
    AUTHOR = {Petrovi\'{c}, Vojislav},
     TITLE = {Antidirected {H}amiltonian circuits in tournaments},
 BOOKTITLE = {Graph theory ({N}ovi {S}ad, 1983)},
     PAGES = {259--269},
 PUBLISHER = {Univ. Novi Sad, Novi Sad},
      YEAR = {1984},
   MRCLASS = {05C45 (05C20)},
  MRNUMBER = {751454},
MRREVIEWER = {Carsten Thomassen},
}

@article {thomason1986paths,
    AUTHOR = {Thomason, Andrew},
     TITLE = {Paths and cycles in tournaments},
   JOURNAL = {Trans. Amer. Math. Soc.},
  FJOURNAL = {Transactions of the American Mathematical Society},
    VOLUME = {296},
      YEAR = {1986},
    NUMBER = {1},
     PAGES = {167--180},
      ISSN = {0002-9947},
   MRCLASS = {05C20 (05C38)},
  MRNUMBER = {837805},
MRREVIEWER = {Brian Alspach},
       DOI = {10.2307/2000567},
       URL = {https://www.ams.org/journals/tran/1986-296-01/S0002-9947-1986-0837805-6/},
}

@article {havet2000oriented,
    AUTHOR = {Havet, Fr\'{e}d\'{e}ric and Thomass\'{e}, St\'{e}phan},
     TITLE = {Oriented {H}amiltonian paths in tournaments: a proof of
              {R}osenfeld's conjecture},
   JOURNAL = {J. Combin. Theory Ser. B},
  FJOURNAL = {Journal of Combinatorial Theory. Series B},
    VOLUME = {78},
      YEAR = {2000},
    NUMBER = {2},
     PAGES = {243--273},
      ISSN = {0095-8956},
   MRCLASS = {05C45 (05C20)},
  MRNUMBER = {1750898},
MRREVIEWER = {Moshe Rosenfeld},
       DOI = {10.1006/jctb.1999.1945},
       URL = {https://www.sciencedirect.com/science/article/pii/S0095895699919457},
}

@article {havet2000orientedcycles,
    AUTHOR = {Havet, Fr\'{e}d\'{e}ric},
     TITLE = {Oriented {H}amiltonian cycles in tournaments},
   JOURNAL = {J. Combin. Theory Ser. B},
  FJOURNAL = {Journal of Combinatorial Theory. Series B},
    VOLUME = {80},
      YEAR = {2000},
    NUMBER = {1},
     PAGES = {1--31},
      ISSN = {0095-8956},
   MRCLASS = {05C20 (05C38 05C45)},
  MRNUMBER = {1778196},
MRREVIEWER = {Brenda J. Latka},
       DOI = {10.1006/jctb.2000.1959},
       URL = {https://www.sciencedirect.com/science/article/pii/S0095895600919592},
}

@online{zein2022oriented,
      title={Oriented {H}amiltonian cycles in tournaments: a proof of {R}osenfeld's conjecture},
      author={Ayman El Zein},
      year={2022},
      eprint={2204.11211},
      archivePrefix={arXiv},
      primaryClass={math.CO}
}

@article {hell2002antidirected,
    AUTHOR = {Hell, P. and Rosenfeld, M.},
     TITLE = {Antidirected {H}amiltonian paths between specified vertices of
              a tournament},
   JOURNAL = {Discrete Appl. Math.},
  FJOURNAL = {Discrete Applied Mathematics. The Journal of Combinatorial
              Algorithms, Informatics and Computational Sciences},
    VOLUME = {117},
      YEAR = {2002},
    NUMBER = {1-3},
     PAGES = {87--98},
      ISSN = {0166-218X},
   MRCLASS = {05C45 (05C20 05C85)},
  MRNUMBER = {1881270},
       DOI = {10.1016/S0166-218X(01)00197-4},
       URL = {https://www.sciencedirect.com/science/article/pii/S0166218X01001974},
}

@article {montgomery2022transversal,
    AUTHOR = {Montgomery, Richard and M\"{u}yesser, Alp and Pehova, Yani},
     TITLE = {Transversal factors and spanning trees},
   JOURNAL = {Adv. Comb.},
  FJOURNAL = {Advances in Combinatorics},
      YEAR = {2022},
     PAGES = {Paper No. 3, 25},
   MRCLASS = {05C35},
  MRNUMBER = {4451150},
       DOI = {10.19086/aic.2022.3},
       URL = {https://www.advancesincombinatorics.com/article/35484-transversal-factors-and-spanning-trees},
}

@article {havet2000median,
    AUTHOR = {Havet, Fr\'{e}d\'{e}ric and Thomass\'{e}, St\'{e}phan},
     TITLE = {Median orders of tournaments: a tool for the second
              neighborhood problem and {S}umner's conjecture},
   JOURNAL = {J. Graph Theory},
  FJOURNAL = {Journal of Graph Theory},
    VOLUME = {35},
      YEAR = {2000},
    NUMBER = {4},
     PAGES = {244--256},
      ISSN = {0364-9024},
   MRCLASS = {05C20 (05C05 05C10)},
  MRNUMBER = {1791347},
MRREVIEWER = {Brenda J. Latka},
       DOI = {10.1002/1097-0118(200012)35:4<244::AID-JGT2>3.0.CO;2-H},
       URL = {https://onlinelibrary.wiley.com/doi/abs/10.1002/1097-0118(200012)35:4<244::AID-JGT2>3.0.CO;2-H},
}

@article {aharoni2020rainbow,
    AUTHOR = {Aharoni, Ron and DeVos, Matt and Gonz\'{a}lez Hermosillo de la Maza, Sebasti\'{a}n and Montejano, Amanda and \v{S}\'{a}mal, Robert},
     TITLE = {A rainbow version of {M}antel's theorem},
   JOURNAL = {Adv. Comb.},
  FJOURNAL = {Advances in Combinatorics},
      YEAR = {2020},
     PAGES = {Paper No. 2, 12},
   MRCLASS = {05C35},
  MRNUMBER = {4125343},
MRREVIEWER = {Andrzej Grzesik},
       DOI = {10.19086/aic.12043},
       URL = {https://www.advancesincombinatorics.com/article/12043-a-rainbow-version-of-mantel-s-theorem},
}

@article {barany1982generalization,
    AUTHOR = {B\'{a}r\'{a}ny, Imre},
     TITLE = {A generalization of {C}arath\'{e}odory's theorem},
   JOURNAL = {Discrete Math.},
  FJOURNAL = {Discrete Mathematics},
    VOLUME = {40},
      YEAR = {1982},
    NUMBER = {2-3},
     PAGES = {141--152},
      ISSN = {0012-365X},
   MRCLASS = {52A35},
  MRNUMBER = {676720},
MRREVIEWER = {Marilyn Breen},
       DOI = {10.1016/0012-365X(82)90115-7},
       URL = {https://www.sciencedirect.com/science/article/pii/0012365X82901157},
}

@online {kalai2009colorful,
    AUTHOR = {Kalai, Gil},
     TITLE = {Colorful Caratheodory revisited},
      YEAR = {2009},
       URL = {https://gilkalai.wordpress.com/2009/03/15/colorful-caratheodory-revisited/}
}

@article {joos2020rainbow,
    AUTHOR = {Joos, Felix and Kim, Jaehoon},
     TITLE = {On a rainbow version of {D}irac's theorem},
   JOURNAL = {Bull. Lond. Math. Soc.},
  FJOURNAL = {Bulletin of the London Mathematical Society},
    VOLUME = {52},
      YEAR = {2020},
    NUMBER = {3},
     PAGES = {498--504},
      ISSN = {0024-6093},
   MRCLASS = {05C70 (05C38 05C45)},
  MRNUMBER = {4171383},
MRREVIEWER = {Vahan V. Mkrtchyan},
       DOI = {10.1112/blms.12343},
       URL = {https://londmathsoc.onlinelibrary.wiley.com/doi/abs/10.1112/blms.12343},
}

@article {kalai2005topological,
    AUTHOR = {Kalai, Gil and Meshulam, Roy},
     TITLE = {A topological colorful {H}elly theorem},
   JOURNAL = {Adv. Math.},
  FJOURNAL = {Advances in Mathematics},
    VOLUME = {191},
      YEAR = {2005},
    NUMBER = {2},
     PAGES = {305--311},
      ISSN = {0001-8708},
   MRCLASS = {52A35 (55N35)},
  MRNUMBER = {2103215},
MRREVIEWER = {Yu. A. Shashkin},
       DOI = {10.1016/j.aim.2004.03.009},
       URL = {https://www.sciencedirect.com/science/article/pii/S0001870804000969},
}

@article {cheng2023rainbow,
    AUTHOR = {Cheng, Yangyang and Han, Jie and Wang, Bin and Wang, Guanghui},
     TITLE = {Rainbow spanning structures in graph and hypergraph systems},
   JOURNAL = {Forum Math. Sigma},
  FJOURNAL = {Forum of Mathematics. Sigma},
    VOLUME = {11},
      YEAR = {2023},
     PAGES = {Paper No. e95, 20},
      ISSN = {2050-5094},
   MRCLASS = {05C35 (05C65)},
  MRNUMBER = {4657062},
       DOI = {10.1017/fms.2023.92},
       URL = {https://doi.org/10.1017/fms.2023.92},
}

@article {aharoni2017rainbow,
    AUTHOR = {Aharoni, Ron and Howard, David},
     TITLE = {A rainbow {$r$}-partite version of the {E}rd\H{o}s-{K}o-{R}ado
              theorem},
   JOURNAL = {Combin. Probab. Comput.},
  FJOURNAL = {Combinatorics, Probability and Computing},
    VOLUME = {26},
      YEAR = {2017},
    NUMBER = {3},
     PAGES = {321--337},
      ISSN = {0963-5483},
   MRCLASS = {05D05},
  MRNUMBER = {3628907},
MRREVIEWER = {Rui Gu},
       DOI = {10.1017/S0963548316000353},
       URL = {https://www.cambridge.org/core/journals/combinatorics-probability-and-computing/article/abs/rainbow-rpartite-version-of-the-erdoskorado-theorem/645E07E1715A0CA2E98E64D039E1451C},
}

@online {chakraborti2023bandwidth,
    AUTHOR = {Chakraborti, Debsoumya and Im, Seonghyuk and Kim, Jaehoon and Liu, Hong},
     TITLE = {A bandwidth theorem for graph transversals},
      YEAR = {2023},
    EPRINT = {2302.09637},
    archivePrefix = {arXiv},
    primaryClass = {math.CO},
}

@online {cheng2021transversal,
    author = {Cheng, Yangyang and Han, Jie and Wang, Bin and Wang, Guanghui and Yang, Donglei},
    title = {Transversal {H}amilton cycle in hypergraph systems},
    year = {2021},
    EPRINT = {2111.07079},
    archivePrefix = {arXiv},
    primaryClass = {math.CO},
}

@article {huang1994relations,
    AUTHOR = {Huang, Rosa and Rota, Gian-Carlo},
     TITLE = {On the relations of various conjectures on {L}atin squares and
              straightening coefficients},
   JOURNAL = {Discrete Math.},
  FJOURNAL = {Discrete Mathematics},
    VOLUME = {128},
      YEAR = {1994},
    NUMBER = {1-3},
     PAGES = {225--236},
       DOI = {10.1016/0012-365X(94)90114-7},
       URL = {https://www.sciencedirect.com/science/article/pii/0012365X94901147},
}

@online {pokrovskiy2020rota,
    author = {Pokrovskiy, Alexey},
    title = {Rota's {B}asis {C}onjecture holds asymptotically},
    year = {2020},
    EPRINT = {2008.06045},
    archivePrefix = {arXiv},
    primaryClass = {math.CO},
  }

@article {cheng2021pancyclicity,
    AUTHOR = {Cheng, Yangyang and Wang, Guanghui and Zhao, Yi},
     TITLE = {Rainbow pancyclicity in graph systems},
   JOURNAL = {Electron. J. Combin.},
  FJOURNAL = {Electronic Journal of Combinatorics},
    VOLUME = {28},
      YEAR = {2021},
    NUMBER = {3},
     PAGES = {Paper No. 3.24, 9},
       DOI = {10.37236/9033},
       URL = {https://www.combinatorics.org/ojs/index.php/eljc/article/view/v28i3p24},
}

@online {chakraborti2023hamilton,
    AUTHOR = {Chakraborti, Debsoumya and Kim, Jaehoon and Lee, Hyunwoo and Seo, Jaehyeon},
     TITLE = {Hamilton transversals in tournaments},
    Note = {accepted by Combinatorica},
      YEAR = {2023},
    EPRINT = {2307.00912},
    archivePrefix = {arXiv},
    primaryClass = {math.CO},
}

@article {gupta2023general,
    AUTHOR = {Gupta, Pranshu and Hamann, Fabian and M\"{u}yesser, Alp and
              Parczyk, Olaf and Sgueglia, Amedeo},
     TITLE = {A general approach to transversal versions of {D}irac-type
              theorems},
   JOURNAL = {Bull. Lond. Math. Soc.},
  FJOURNAL = {Bulletin of the London Mathematical Society},
    VOLUME = {55},
      YEAR = {2023},
    NUMBER = {6},
     PAGES = {2817--2839},
      ISSN = {0024-6093,1469-2120},
   MRCLASS = {05D40 (05D15)},
  MRNUMBER = {4689555},
       DOI = {10.1112/blms.12896},
       URL = {https://doi.org/10.1112/blms.12896},
}

@online {tang2023rainbow,
    AUTHOR = {Tang, Yucong and Wang, Bin and Wang, Guanghui and Yan, Guiying},
     TITLE = {Rainbow {H}amilton cycle in hypergraph system},
      YEAR = {2023},
    EPRINT = {2302.00080},
    archivePrefix = {arXiv},
    primaryClass = {math.CO},
}

@online {ferber2022diractype,
    AUTHOR = {Ferber, Asaf and Han, Jie and Mao, Dingjia},
     TITLE = {Dirac-type Problem of Rainbow matchings and Hamilton cycles in Random Graphs},
      YEAR = {2022},
    EPRINT = {2211.05477},
    archivePrefix = {arXiv},
    primaryClass = {math.CO},
}

@article {kuhn2011sumner,
    AUTHOR = {K\"{u}hn, Daniela and Mycroft, Richard and Osthus, Deryk},
     TITLE = {A proof of {S}umner's universal tournament conjecture for
              large tournaments},
   JOURNAL = {Proc. Lond. Math. Soc. (3)},
  FJOURNAL = {Proceedings of the London Mathematical Society. Third Series},
    VOLUME = {102},
      YEAR = {2011},
    NUMBER = {4},
     PAGES = {731--766},
      ISSN = {0024-6115,1460-244X},
   MRCLASS = {05C20 (05C05 05C35)},
  MRNUMBER = {2793448},
MRREVIEWER = {K.\ B.\ Reid},
       DOI = {10.1112/plms/pdq035},
       URL = {https://doi-org-ssl.access.yonsei.ac.kr/10.1112/plms/pdq035},
}

@online {anastos2023robust,
    AUTHOR = {Anastos, Michael and Chakraborti, Debsoumya},
     TITLE = {Robust Hamiltonicity in families of Dirac graphs},
      YEAR = {2023},
    EPRINT = {2309.12607},
    archivePrefix = {arXiv},
    primaryClass = {math.CO},
}

@online {babinski2023directed,
    AUTHOR = {Babi{\'n}ski, Sebastian and Grzesik, Andrzej and Prorok, Magdalena},
     TITLE = {Directed graphs without rainbow triangles},
      YEAR = {2023},
    EPRINT = {2308.01461},
    archivePrefix = {arXiv},
    primaryClass = {math.CO},
}

@online {gerbner2024directed,
    AUTHOR = {Gerbner, Daniel and Grzesik, Andrzej and Palmer, Cory and Prorok, Magdalena},
     TITLE = {Directed graphs without rainbow stars},
      YEAR = {2024},
    EPRINT = {2402.01028},
    archivePrefix = {arXiv},
    primaryClass = {math.CO},
}

@article {bang-jensen1992polynomial,
    AUTHOR = {Bang-Jensen, J\o rgen and Manoussakis, Yannis and Thomassen,
              Carsten},
     TITLE = {A polynomial algorithm for {H}amiltonian-connectedness in
              semicomplete digraphs},
   JOURNAL = {J. Algorithms},
  FJOURNAL = {Journal of Algorithms. Cognition, Informatics and Logic},
    VOLUME = {13},
      YEAR = {1992},
    NUMBER = {1},
     PAGES = {114--127},
      ISSN = {0196-6774},
   MRCLASS = {05C85 (05C45 68R10)},
  MRNUMBER = {1146335},
MRREVIEWER = {Zden\v{e}k\ Ryj\'{a}\v{c}ek},
       DOI = {10.1016/0196-6774(92)90008-Z},
       URL = {https://doi-org-ssl.access.yonsei.ac.kr/10.1016/0196-6774(92)90008-Z},
}

@article {hell1983complexity,
    AUTHOR = {Hell, Pavol and Rosenfeld, Moshe},
     TITLE = {The complexity of finding generalized paths in tournaments},
   JOURNAL = {J. Algorithms},
  FJOURNAL = {Journal of Algorithms},
    VOLUME = {4},
      YEAR = {1983},
    NUMBER = {4},
     PAGES = {303--309},
      ISSN = {0196-6774},
   MRCLASS = {05C20 (68Q25 68R10)},
  MRNUMBER = {729226},
       DOI = {10.1016/0196-6774(83)90011-1},
       URL = {https://doi-org-ssl.access.yonsei.ac.kr/10.1016/0196-6774(83)90011-1},
}

\end{document}